\documentclass[12pt]{article}

\usepackage{a4wide}
\usepackage{amssymb}
\usepackage{amsfonts}
\usepackage{amsmath}
\input xy
\xyoption{arrow} \xyoption{matrix}

\date{}

\newtheorem{proposition}{Proposition}[section]
\newtheorem{theorem}[proposition]{Theorem}
\newtheorem{lemma}[proposition]{Lemma}

\newtheorem{corollary}[proposition]{Corollary}

\def\der{\partial }

\def\nFM0{{\nu }_{F,M_0}}
\def\nFN0{{\nu }_{F,N_0}}
\def\nGN0{{\nu }_{G,N_0}}

\def\N0{ {\bf N}_0 }

\def\g{\gamma}

\def\ra{\rightarrow}

\def\Xpm{X^{\pm }}

\def\s{\sigma}
\def\Z{\mathbb{Z}}

\def\l1{{\lambda}_1}

\def\a{\alpha}
\def\a0{ {\alpha }_0}
\def\a1{ {\alpha }_1}

\def\l{\lambda}
\def\o{\omega}
\def\nFGM0{{\nu }_{F,G,M_0}}

\def\nFN0{{\nu}_{F,N_0}}

\def\sm{{\sigma}^m}

\def\sm1{{\sigma}^{-1}}

\def\smtp1{{\sigma}^{-t+1}}

\def\o{\omega }
\def\S1{S^{-1}}

\def\Xpm1{X^{\pm 1}_1}

\def\sPM1{{\sigma }^{\pm 1}}
\def\sMP1{{\sigma }^{\mp 1 }}

\def\di{{\rm d.ind}}

\def\L{\Lambda}

\def\G{\Gamma}

\def\OO{{\cal O}}

\def\CD{{\cal D}}

\def\Ytm1{Y^{t-1}}
\def\Yim1{Y^{i-1}}

\def\CM{{\cal M}}
\def\CN{{\cal N}}

\def\CZ{{\cal Z}}

\def\dim{{\rm dim }}

\def\ker{ {\rm ker } }

\def\gcd{ {\rm gcd } }
\def\D{ \Delta }
\def\SL2Z{ {\rm SL}_2({\bf Z}) }

\def\CZ{ {\cal Z}}
\def\CR{ {\cal R}}

\def\Gp1{ G^{1 , 1 } }
\def\P11{ P^{-1 , 1 } }
\def\Pp1{ P^{1 , 1 } }

\def\lcm{{\rm lcm }}

\def\CE{{\cal E}}

\def\nCLsr{{}^\nu\kern-2pt {\cal L}^{\sigma , \rho  }}
\def\nP{{}^\nu \kern-2pt P}
\def\nL{{}^\nu\kern-2pt L}
\def\nLL{{}^\nu\kern-2pt \Lambda}
\def\nPsr{{}^\nu\kern-2pt P^{\sigma , \rho  }}
\def\nLsr{{}^\nu\kern-2pt L^{\sigma , \rho  }}
\def\nuCL{{}^\nu\kern-2pt  {\cal L}}
\def\nCLsr{{}^\nu\kern-2pt {\cal L}^{\sigma , \rho  }}
\def\nCL1m{{}^\nu\kern-2pt {\cal L}^{-1 , 1  }}
\def\x1nu{x^\frac{1}{\nu}}
\def\xm1nu{x^{-\frac{1}{\nu}}}

\def\trdeg{{\rm tr.deg}}

\def\CR{ {\cal R}}

\def\CN{{\cal N}}
\def\ra{\rightarrow }

\def\lcm{{\rm lcm}}

\def\CC{ {\cal C}}
\def\CE{ {\cal E} }

\def\CP{ {\cal P}}

\def\nAM0{{\nu }_{{\cal A},M_0}}
\def\nAN0{{\nu }_{{\cal A},N_0}}

\def\CR{ {\cal R }}
\def\CP{ {\cal P }}

\def\gn{\mathfrak{n}}

\def\SL{{\rm SL}}

\def\di!{\frac{\der^i}{i!}}
\def\dik!{\frac{\der^k_i}{k!}}
\def\Max{{\rm Max}}

\def\N{\mathbb{N}}

\def\0{\overline{0}}
\def\1{\overline{1}}

\def\Ln1{\L_{n,\overline{1}}}

\def\a1{a_{\overline{1}}}

\def\S{\Sigma}

\def\Dec{{\rm Dec}}

\def\vn1{\overrightarrow{n-1}}
\def\CQ{{\cal Q}}

\def\Q{\mathbb{Q}}

\def\C{\mathbb{C}}
\def\Dec{{\rm Dec}}
\def\IrrA{{\rm Irr}(A)}
\def\IrrK[x]{{\rm Irr}(K[x])}
\def\Adm{{\rm Adm}}
\def\df{{\rm def}}
\def\ind{{\rm ind}}

\begin{document}

\author{V. V. \  Bavula %(genRitt.tex)
}

\title{Generalizations of two theorems of Ritt on decompositions of
polynomial maps}

\maketitle

$$Dedicated\; to\; F.\; van\;  Oystaeyen\;on \; the\; occasion\; of\; his \;
60'th\; birthday $$

\begin{abstract}
In 1922, J. F. Ritt \cite{RittTrAMS1922} proved two remarkable
theorems on decompositions of polynomial maps of $\C [x]$ into
irreducible polynomials (with respect to the composition $\circ$ of
maps). Briefly, the first theorem states that in any two
decompositions of a given polynomial function into irreducible
polynomials the number of the irreducible polynomials and their
degrees are the same (up to order). The second theorem gives four
types of transformations of how to obtain all the decompositions
from a given one. In 1941, H. T. Engstrom \cite{EngstromAJM1941}
and, in 1942, H. Levi \cite{LeviAJM1942} generalized respectively
the first and the second theorem to polynomial maps over an
arbitrary field $K$ of characteristic
zero. %The proof of H. T. Engstrom is not complete.
The aim of the paper is
%to bridge the gap and
 to generalize the two theorems of J. F. Ritt to a more general situation: for, so-called,
{\em reduction} monoids   ($(K[x], \circ )$ and $(K[x^2]x, \circ)$
are examples of reduction monoids). In particular, analogues of the
two theorems of J. F. Ritt hold for the monoid $(K[x^2]x, \circ)$ of
odd polynomials. It is shown that, in general, the two theorems of
J. F. Ritt fail for the cusp $(K+K[x]x^2, \circ )$ but their
analogues are still true for decompositions of {\em maximal length}
of {\em regular} elements of the cusp.

{\em Key Words:  the two theorems of Ritt, Ritt transformations,
composition of polynomial maps, cusp transformations, irreducible
map, the length and defect of a polynomial. }

 {\em Mathematics subject classification
2000:  12F20, 14H37, 14R10.}

$${\bf Contents}$$
\begin{enumerate}
\item Introduction. \item Generalizations of the two theorems of
J. F. Ritt.\item Analogues of the two theorems of J. F. Ritt for
the  cusp.
\end{enumerate}
\end{abstract}

%%%%%%%%%%%%%%%%%% SECTION 1 %%%%%%%%%%%%%%%%%%%%%%%%

\section{Introduction}
In this paper, $K$ is a field of characteristic zero and $K[x]$ is
a polynomial algebra over the field $K$ in a single variable $x$.
The polynomial algebra $K[x]$ is a monoid, $(K[x], \circ )$, where
$\circ$ is the {\em composition} of polynomial functions, $(a\circ
b)(x):= a(b(x))$, and $x$ is the identity element of the monoid
$K[x]$. An element $u$ of the monoid $K[x]$ is a unit iff $\deg
(u)=1$. The group of units of the monoid $(K[x], \circ )$ is
denoted by $K[x]^*$.

A polynomial $a\in K[x]$ is said to be {\em irreducible} (or {\em
prime} or {\em indecomposable}) if $\deg (a)>1$ and the polynomial
$a$ is not a composition of two non-units, i.e. $a$ is an
irreducible element of the monoid $(K[x], \circ )$. This concept of
irreducibility should not be confused with the concept of
irreducibility of the multiplicative monoid $(K[x], \cdot )$ which
is {\em not} used in the paper. A polynomial which is not
irreducible is said to be {\em reducible} or {\em composite}. When
$K= \C$ composite polynomials were studied by J. F. Ritt
\cite{RittTrAMS1922}. He proved two theorems that completely
describe the decompositions composite polynomials may possess. His
{\em first} theorem states: {\em any two decompositions of a given
polynomial  of $\C [x]$ into irreducible polynomials contain the
same number of polynomials; the degrees of the polynomials in one
decomposition are the same as those in the other, except, perhaps,
for the order in which they occur.}

Two decompositions of a polynomial $a$ into irreducible
polynomials
$$ a = p_1\circ \cdots \circ p_r= q_1\circ \cdots \circ q_r$$
are called {\em equivalent} if there exist $r-1$ polynomials of
the first degree $u_1, \ldots , u_{r-1}$ such that
$$ q_1= p_1\circ u_1, \;\; q_2= u_1^{-1} \circ p_2\circ u_2,
\ldots , q_{r-1}= u^{-1}_{r-2} \circ p_{r-1}\circ u_{r-1}, \;\;
q_r= u^{-1}_{r-1}\circ p_r.$$ Suppose that in a decomposition of
the polynomial $a$ into irreducible polynomials %\marginpar{ap1p1}
\begin{equation}\label{ap1p1}
a= p_1\circ \cdots \circ p_r
\end{equation}
there is an adjacent pair of irreducible polynomials
$$ p_i=\l_1\circ \pi_1\circ \l_2, \;\; p_{i+1} = \l_2^{-1}\circ
\pi_2\circ \l_3$$ where $\l_1$, $\l_2$ and $\l_3$ are polynomials
of degree 1 and where $\pi_1$ and $\pi_2$, of unequal degrees $m$
and $n$, respectively, are of any of the following three types:

\begin{eqnarray*}
 &(a)& \pi_1=T_m, \;\;\;\; \pi_2=T_n, \\
 &(b)& \pi_1=x^m, \;\;\;\; \pi_2=x^rg(x^m), \\
 &(c)& \pi_1=x^rg^n, \;\; \pi_2=x^n,
\end{eqnarray*}
where $g=g(x)$ is a polynomial, $T_n$ is the trigonometric
polynomial, $T_n(\cos \, t) := \cos (nt)$.

Then, for the polynomial $a$ we have a decomposition distinct from
(\ref{ap1p1}),

 %\marginpar{ap1p2}
\begin{equation}\label{ap1p2}
a= p_1\circ \cdots \circ p_{i-1} \circ p_i^*\circ p_{i+1}^*\circ
p_{i+2}\circ \cdots \circ p_r
\end{equation}
where respectively to the three cases above the polynomials
$p_i^*$ and $p_{i+1}^*$ are as follows:

\begin{eqnarray*}
 &(a)& p_i^*=\l_1\circ T_n, \;\;\;\;\;\;\;\;\, p_{i+1}^*=T_m\circ \l_3, \\
 &(b)& p_i^*=\l_1\circ [x^rg^m], \;\; p_{i+1}^*=x^m\circ \l_3, \\
 &(c)& p_i^*=\l_1\circ x^n, \;\;\;\;\;\;\;\;\, p_{i+1}^*=[x^rg(x^n)]\circ \l_3.
\end{eqnarray*}

Clearly, $\deg (p_i^*) = \deg (p_{i+1})= n$ and $\deg (p_{i+1}^* )
= \deg (p_i)= m$.

$\noindent $

The {\em second} theorem of J. F. Ritt states: {\em if $a\in \C
[x]$ has two distinct decompositions into irreducible polynomials,
we can pass from either to a decomposition equivalent to the other
by repeated steps of the three types just indicated.}

$\noindent $

 He writes in his paper, p. 53: {\em ``The analogous problem for
 fractional rational functions  is much more difficult. There is a
 much greater variety of possibilities, as one sees, without going
 far, on considering the formulas for the transformation of the
 periods of the elliptic functions. There are even cases in which
 the number of prime functions in one decomposition is different
 from that in another.''} We will see later in the paper that the
 situation is similar for the cusp.

J. F.  Ritt's approach is based on the monodromy group associated
with the equation $f(x)-y=0$.

Later H. T. Engstrom \cite{EngstromAJM1941} and H. Levi
\cite{LeviAJM1942} proved respectively the first and the second
theorem of J. F. Ritt for the polynomial algebra $K[x]$ where $K$ is
a field of characteristic zero.  Their methods are algebraic.
%In the proof of the generalization of the second theorem of J. F.
%Ritt, H. Levi used the generalization of H. T. Engstrom (Theorem
%4.1, \cite{EngstromAJM1941}) of the first theorem of J. F. Ritt.
%The proof of H. T. Engstrom of Theorem 4.1,
%\cite{EngstromAJM1941}, is not complete: the proof is valid if
%$|C|>1$ (in the notation of Theorem 4.1, \cite{EngstromAJM1941});
%if $|C|=1$ his induction does not work. In fact, his arguments are
%just a tautology (see the last two sentences of his proof), and
%there are no arguments to cover (or to exclude) this case. The
%case $|C|=1$ may well occur in the proof. Example: $C=1$,
%$A=g_s=T_p$ and $f_r=B=T_q$ where $T_p$ and $T_q$ are the
%trigonometric polynomials when $p$ and $q$ are distinct prime
%numbers.
% Another example is when $C=1$, $A=g_s=x^p$ and $f_r= B=x^q$ where
% $p$ and $q$ are distinct prime numbers. There are more
% sophisticated examples of this sort.

It is known that the theorems of J. F. Ritt are false in prime
characteristic \cite{Dorey-Whaples1974},
\cite{Gutierrez-Sevilla2006}, but the first theorem is true for,
so-called, tame polynomials \cite{Fried-MacRae1969},
\cite{Schinzel1982}. For some generalizations, applications and
connections with the two theorems of J. F. Ritt the reader is
referred to \cite{Avanzi-Zannier2003, Bilu-Tichy2000, Binder95,
Eigenthalter-Woracek1995, Gutierrez-Sevilla2006, Pakovicharxiv2007,
Schinzel1982, Tortrat1988, Zannier1993, Zieve2007}.

The goal of this paper is to generalize the two theorems of J. F.
Ritt to a more general situation (for, so-called, {\em reduction}
monoids - see Section \ref{G2THMR} for a definition; $(K[x], \circ
)$ and  $(K[x^2]x, \circ )$ are  reduction monoids).
% and, as a
%consequence, to bridge the mentioned gap.
 The advantage of our
method is that generalizations of the two theorems are proved in one
go.

 For a natural number $r$, let $S_r$ be the symmetric
group. For reduction monoids (the definition is given in Section
\ref{G2THMR}), the first and the second statement of the following
theorem are generalizations of the first and the second theorem of
J. F. Ritt, respectively.  The first statement is precisely the same
as the first theorem of J. F. Ritt, but the second statement
contains only `half' of the second theorem of J. F. Ritt,  as the
second part of the second theorem of J. F. Ritt classifies all the
situations $p_ip_{i+1}=p_i'p_{i+1}'$ for the monoid $(\C [x], \circ
)$.

\begin{theorem}\label{b31Aug7}%\marginpar{b31Aug7}
Let $\CM$ be a reduction monoid, $\CM^*$ be its group of units,
$a\in \CM$ with $|a|>1$, and $a=p_1\cdots p_r= q_1\cdots q_s$ be
two decompositions of the element $a$ into irreducible factors.
Then
\begin{enumerate}
\item $r=s$ and $|p_1|=|q_{\s (1)}|, \ldots , |p_r|=|q_{\s (r)}|$
for a permutation $\s \in S_r$; and \item if the decompositions
are distinct then one can be obtained from the other by finitely
many transformations on adjacent irreducible factors of the
following two types:
\begin{enumerate}
 \item $p_1\cdots p_ip_{i+1}\cdots  p_r=p_1\cdots (p_iu)(u^{-1}p_{i+1})\cdots
 p_r$ where $u\in \CM^*$,
\item $p_1\cdots p_ip_{i+1}\cdots  p_r=p_1\cdots
p_i'p_{i+1}'\cdots p_r$ where $p_ip_{i+1}=p_i'p_{i+1}'$, the
numbers $|p_i|$ and $|p_{i+1}|$ are co-prime, $|p_i|=|p_{i+1}'|$
and $|p_{i+1}|=|p_i'|$.
\end{enumerate}
\end{enumerate}
\end{theorem}

Consider the submonoid $(\OO := K[x^2]x, \circ)$ of odd polynomials
of the monoid $(K[x], \circ )$.

\begin{theorem}\label{2Nov7}%\marginpar{2Nov7}
Let $K$ be a field of
characteristic zero. Then the monoid $\OO $ is a reduction monoid
where $|\cdot |= \deg $.
\end{theorem}

The group $\OO^*$ of units of the monoid $\OO$ is equal to the group
$\{ \l x  \, | \, \l \in K^*\}$ where  $K^*:=K\backslash \{ 0\}$.
The first two statements of the next corollary follow at once from
Theorems \ref{b31Aug7} and \ref{2Nov7}; statement 3 follows from the
second theorem of J. F. Ritt but not in a straightforward way as
many additional results are used in its proof: Theorem \ref{8Nov7},
Lemma \ref{3Nov7}, Lemma \ref{c8Nov7} (see Section \ref{G2THMR} for
detail).

\begin{corollary}\label{a2Nov7}%\marginpar{a2Nov7}
Let $K$ be a field of characteristic zero, $a\in \OO$ with $\deg
(a)>1$, and $a=p_1\circ \cdots \circ p_r= q_1\circ \cdots \circ q_s$
be two decompositions of the element $a$ into irreducible factors of
   the monoid $\OO$. Then
\begin{enumerate}
\item $r=s$ and $\deg (p_1)=\deg (q_{\s (1)}), \ldots , \deg (p_r)=\deg (q_{\s (r)})$
for a permutation $\s \in S_r$; and \item if the decompositions are
distinct then one can be obtained from the other by finitely many
transformations on adjacent irreducible factors of the following two
types:
\begin{enumerate}
 \item $p_1\circ \cdots \circ p_i\circ p_{i+1}\circ \cdots  \circ p_r=p_1\circ \cdots \circ (p_i\circ u)\circ
  (u^{-1}\circ p_{i+1})\circ \cdots \circ
 p_r$ where $u\in \OO^*$,
\item $p_1\circ \cdots \circ p_i\circ p_{i+1}\circ \cdots \circ  p_r=p_1\circ \cdots\circ
p_i^*\circ p_{i+1}^*\circ \cdots \circ p_r$ where $$p_i\circ
p_{i+1}=p_i^*\circ p_{i+1}^*,$$ the degrees $\deg (p_i)$ and $\deg
(p_{i+1})$ are co-prime, $\deg (p_i)=\deg (p_{i+1}^*)$ and $\deg
(p_{i+1})=\deg (p_i^*)$.
\end{enumerate}
\item There are only the following options for the pairs $P:= (p_i,
p_{i+1})$ and $P^*:= (p_i^*, p_{i+1}^*)$:
\begin{enumerate}
\item $P=(T_n, T_m)$ and $P^*=(T_m, T_n)$ where $n$ and $m$ are odd
distinct primes,
\item $P=(x^t[\alpha (x^2)]^s, x^s)$ and $P^*=(x^s, x^t\alpha (x^{2s}))$,
\item $P=(x^s, x^t\alpha (x^{2s}))$ and $P^*=(x^t[\alpha (x^2)]^s, x^s)$,
\end{enumerate}
where $s$ is an odd prime number, $t$ is an odd number,  and $\alpha
\in K[x]\backslash K$ with $\alpha (0)\neq 0$.
\end{enumerate}
\end{corollary}

Up to my knowledge, the monoid $\OO$ is the only example distinct
from $K[x]$ for which  (analogues of) the two theorems of J. F. Ritt
hold. It would be interesting to find more examples (the definition
of reduction monoid is very arithmetical). It is a curious fact that
the monoid $\OO$, in fact,  comes  from non-commutative situation.
The monoid $\OO$ is the monoid of all central algebra endomorphisms
of a certain localization of the quantum plane which is a {\em
non-commutative} algebra (see Section \ref{G2THMR} for detail). It
would be interesting to find more reduction monoids coming from
non-commutative situation (and as a result to obtain analogues of
the two theorems of J. F. Ritt for them). Notice that in the
definition of reduction monoid $\CM$ is not necessarily a {\em
commutative} algebra, it is just an abelian group. Moreover, in the
case of the odd polynomials,  $\OO$ is  not even an algebra.

The cusp submonoid $(K+K[x]x^2, \circ )$ of $(K[x], \circ )$ looks
similar to the monoid $\OO$ but for it situation is completely
different. In particular, the cusp submonoid is not a reduction
monoid.

Till the end of this section let $K$ be an {\em algebraically
closed} field of characteristic zero and let $A$ be the subalgebra
of the polynomial algebra $K[x]$ generated by the monomials $x^2$
and $x^3$. The algebra $A=K+K[x]x^2$ is isomorphic to the algebra of
regular functions on the cusp $s^2= t^3$. It is obvious that $(A,
\circ )$ is a sub-semi-group of $(K[x], \circ )$. For a polynomial
$a\in K[x]$ of degree $\deg (a)>1$, let $\Dec (a)$ be the set of all
decompositions of the polynomial $a$ into irreducible polynomials
 of $K[x]$ (with respect to $\circ$). The {\em length} $l(a)$ of the polynomial $a\in K[x]$ is
the number of irreducible polynomials in any decomposition of $\Dec
(a)$. Similarly, for a polynomial $a\in A\backslash K$, let
$\Dec_A(a)$  be the set of all decompositions of the polynomial $a$
into irreducible polynomials of $A$. The natural number
$$ l_A(a):= \max \{ r\, | \, p_1\circ \cdots \circ p_r\in
\Dec_A(a)\}$$ is called the $A$-{\em length} of the element $a$.
It is obvious that
$$l_A(a) \leq l(a).$$
In general,  this inequality is strict (Corollary \ref{10Sep7}). An
element $a\in A$ is called {\em regular} (respect. {\em irregular})
if $l_A(a) = l(a)$ (resp. $l_A(a) < l(a)$). The are plenty of
elements of both types. Moreover, if $a$ is irregular then $a\circ
(x+\l )$ is regular for some $\l\in K$.  A decomposition
$$p_1\circ \cdots \circ p_{l_A(a)}\in
\Dec_A(a)$$ is called a {\em decomposition of maximal length} or a
{\em maximal decomposition} for the element $a$. Let $\Max (a)$ be
the set of all maximal decompositions for $a$. Clearly, $\Max (a)
\subseteq \Dec_A (a)$, but, in general, $\Max (a) \neq \Dec_A (a)$,
see (\ref{MaDe}). Lemma \ref{22Sep7} describes the set $\Max (a)$.

In general, the number of irreducible polynomials in decomposition
into irreducible polynomials of an element of $A$ is non-unique
(Lemma \ref{b9Sep7}); moreover, it can vary  greatly. So, for the
cusp the two theorems of J. F. Ritt do not hold. Therefore, the cusp
is not a reduction monoid. Nevertheless, for decompositions of
maximal length of each regular element $a$ of $A$ analogues of the
two theorems do hold -- Theorem \ref{15Sep7} and Theorem
\ref{A15Sep7} if $K$ is algebraically closed (if $K$ is not
algebraically closed then, in general, Theorem \ref{A15Sep7} does
not hold).

\begin{theorem}\label{15Sep7}%\marginpar{15Sep7}
Let $K$ be a field of characteristic zero, $a$ be a regular element
of $A$ such that $a\not\in K$, and $$a= p_1\circ \cdots \circ p_r=
q_1\circ \cdots \circ q_r$$ be two decompositions of maximal length
of the element $a$ into irreducible polynomials of $A$. Then
$$\deg (p_1)=\deg (q_{\s (1)}), \ldots , \deg (p_r)=\deg (q_{\s (r)})$$ for a
permutation $\s \in S_r$.
\end{theorem}

Theorem \ref{15Sep7} follows from the first theorem  of J. F. Ritt
(or from Theorem  \ref{A15Sep7}).  In general, for irregular
elements Theorem \ref{15Sep7} is not true (Proposition
\ref{21Sep7}), i.e. the invariance of degrees (up to permutation)
does not hold. The next theorem is  an analogue of the second
theorem of J. F. Ritt for regular elements.
 A new moment is that the transformations (Adm), ($\CC a$), ($\CC
b$) and  ($\CC c$) are defined on {\em three} adjacent elements
rather than two as in the second theorem of J. F. Ritt.

\begin{theorem}\label{A15Sep7}%\marginpar{A15Sep7}
Let $K$ be an algebraically closed field of characteristic zero,
$a$ be a regular element of $A$ such that $a\not\in K$,  and $X, Y
\in \Max (a)$. Then the decomposition $Y$ can be obtained from the
decomposition $X$ by finitely many transformations of the
following four types: (Adm), ($\CC a$), ($\CC b$) and  ($\CC c$),
see below.
\end{theorem}

For a non-scalar polynomial $f$  of $K[x]$, a polynomial $\l +\mu
x$ of degree 1 is called an $f$-{\em admissible} polynomial if
$\l$ is a root of the derivative $f':=\frac{df}{dx}$ of $f$.

Let $a\in A\backslash K$ with
 $r:=l_A(a)=l(a)$, and
$Z := p_1\circ \cdots \circ p_i\circ p_{i+1}\circ \cdots \circ
p_r\in \Max (a)$. Consider the following four types of
transformations of the decomposition $Z$ that produce a new
decomposition $Z^*\in \Max (a)$ where
$$ Z^*:=\begin{cases}
p_1\circ \cdots \circ p_{i-1}\circ p_i^*\circ p_{i+1}^*\circ
p_{i+2}^*\circ \cdots \circ p_r& \text{if $i+1<r$},\\
p_1\circ \cdots \circ p_{r-1}^*\circ p_r^*& \text{if $i+1=r$}.\\
\end{cases}
$$

(Adm) In both cases, $p_i^*:= p_i\circ u$ and $ p_{i+1}^*:= u^{-1}
\circ p_{i+1}$ where $u\in K[x]^*$ is  $p_i$-admissible, and
$p_{i+2}^*= p_{i+2}$ if $i+1<r$ ($u^{-1}$ is the inverse of the
element $u$ in the monoid $(K[x], \circ )$, i.e. $u^{-1}$ is the
inverse map  of $u$).

$\noindent $

In the remaining three  cases below,  $\gcd (\deg (p_i), \deg
(p_{i+1}))=1$, all $\l_i\in K[x]^*$, $p$ is a prime number,
polynomials $x^sg^p(x)$ and $x^sg(x^p)$ satisfy the condition that
$g(0)\neq 0$, $\l_i^{-1}$ is the inverse of the element $\l_i$ in
 the monoid $(K[x], \circ )$.

$\noindent $

($\CC a$) If $i+1<r$, $p_i=\l_1\circ T_k\circ \l_2$ and $p_{i+1}=
\l_2^{-1} \circ T_l\circ \l_3$ where $k$ and $l$ are distinct {\em
odd}  prime numbers, $\l_2$ is $T_k$-admissible and $\l_3$ is
$T_l$-admissible, then
$$p_i^*:= \l_1\circ T_l\circ \l_4, \;\; p_{i+1}^*:=
\l_4^{-1} \circ T_k\circ \l_3 \circ \l_5\;\; {\rm and }\;\;
p_{i+2}^*:=\l_5^{-1} \circ p_{i+2}, $$ where $\l_4$ is
$T_l$-admissible and $\l_5$ is $T_k\circ \l_3$-admissible.

$\noindent $

($\CC b$) If  $i+1<r$, $p_i=\l_1\circ x^p$ and
$p_{i+1}=[x^sg(x^p)]\circ \l_2$ where  $\l_2$ is
$x^sg(x^p)$-admissible, then
$$p_i^*:= \l_1\circ
[x^sg^p]\circ \l_3,\;\;  p_{i+1}^*:= \l_3^{-1} \circ x^p\circ
\l_2\circ \l_4\;\; {\rm and }\;\;p_{i+2}^*:=\l_4^{-1}\circ p_{i+2},
$$ where $\l_3$ is $x^sg^p$-admissible and $\l_2\circ \l_4$ is
$x^p$-admissible.

If $i+1=r$, $p_{r-1}=\l_1 \circ x^p$ and $p_r= [x^s g(x^p)]\circ
\l_2$ where $s\geq  2$ and $\l_2\in K^*x$, then
$$p_{r-1}^*:= \l_1 \circ [x^sg^p] \;\; {\rm and }\;\; p_r^*:= x^p\circ
\l_2.$$

$\noindent $

($\CC c$) If  $i+1<r$, $p_i=\l_1\circ [x^sg^p]\circ \l_2$ and
$p_{i+1}=\l_2^{-1}\circ x^p\circ  \l_3$ where $\l_2$ is
$x^sg^p$-admissible and $\l_3$ is $x^p$-admissible,  then
$$p_i^*:= \l_1\circ x^p,
 \;\;  p_{i+1}^*:= [x^sg(x^p)]\circ \l_3\circ \l_4  \;\; {\rm and }\;\; p_{i+2}^*:=\l_4^{-1} \circ p_{i+2}, $$
where $\l_3\circ \l_4$ is $x^s g(x^p)$-admissible.

If $i+1=r$, $p_{r-1}= \l_1\circ x^s g^p$, $s \geq 2$, and $p_r=
x^p\circ \l_2$ where $\l_2$ is $x^p$-admissible, then
$$ p_{r-1}^*:= \l_1\circ x^p\;\; {\rm and}\;\; p_r^*:= [x^s
g(x^p)]\circ \l_2.$$

Decompositions of polynomials with coefficients in a commutative
ring were studied by the author in  \cite{BavulaGMJ2008}.

%%%%%%%%%%%%%%%%%% SECTION 2 %%%%%%%%%%%%%%%%%%%%%%%%

\section{Generalizations of the two theorems of J. F.
Ritt}\label{G2THMR}%\marginpar{G2THMR}

In this section, the two theorems of J. F. Ritt are generalized to a
more general situation. They are proved for reduction monoids
(Theorem \ref{b31Aug7}). The polynomial algebra $K[x]$ is a
reduction monoid with respect to the composition of functions. These
generalizations are inspired by the paper of H. T. Engstrom
\cite{EngstromAJM1941} and we follow some of his ideas. Proofs of
Theorem \ref{b31Aug7}, Theorem \ref{2Nov7} and Corollary
\ref{a2Nov7}.(3) are given.

Natural numbers $i$ and $j$ are called {\em co-prime} (or {\em
relatively prime}) if $\gcd (i,j)=1$.

$\noindent $

{\it Definition}. A multiplicative monoid $\CM$ is called a {\em
reduction} monoid if the following axioms hold for all elements
$a,b,c\in \CM$ (where $\CM^*$ is the group of units of the monoid
$\CM$):

(A1) $\CM$ is a $\Z$-module (i.e. $\CM$ is an abelian group under
$+$) such that $$(a+b)c= ac+bc.$$

(A2) There exists a map $|\cdot |:\CM \ra \N := \{ 0,1, \ldots \}$
such that $$|ab|= |a||b|\;\; {\rm and} \;\; |a+b| \leq \max \{
|a|, |b|\}.$$

(A3) $a\in \CM^*$ iff $|a|=1$.

$\noindent $

(A4) If $ac=bc$ then $a=b$ provided $|c|>1$. % ; and If $ca=cb$
%then $a=b$ provided $|c|>1$.

$\noindent $

(A5) For any elements $a, b\in \CM$ with $|a|>1$ and  $|b|>1$ and,
in addition,  there exists
 an element $x\in \CM a \cap \CM b$ such that $|x|\neq 0$,  there exists an element $c\in \CM$ such that
 $\CM a\cap \CM b = \CM c$ and $|c|=\lcm (|a|, |b|)$.

$\noindent $

(A6) If $\alpha a = \beta b$ with $|\alpha | =i$, $|a|=jk$,
$|\beta |=j$, $|b|=ik$, $ijk\geq 1$,
 and the natural numbers
 $i$ and $j$ are co-prime then $a=a_1c$ and $b=b_1c$ for some elements $a_1$, $b_1$ and $c$ of $\CM$ such that $|c|=k$.

%(A6) There is a set $\{ \th_i \, | \, i\geq 2, \im  \, | \cdot
%|\}$ of elements of $\CM$ such that
% $|\th_i|=i$ and for each element $a\in \CM$ with $|a|=nm$, $n>1$, $m>1$, there exists a {\em unique} element
%$a^*\in \CM$ of type $a^* = u\th_nb$ where $u\in \CM^*$ and
%$|b|=m$ such that $|a-a^*|<|a|-m$. For any two such presentations,
%say $u\th_nb$ and $u'\th_n'b'$: $u\th_n\CM^*= u'\th_n'\CM^*$.

$\noindent $

{\it Example}. $(K[x], \circ )$ is the reduction monoid where
$|\cdot |:= \deg$. The axioms (A1)-(A4) are obvious. The axioms
(A5) and (A6) follow respectively from Theorems 2.2 and 3.1 of the
paper \cite{EngstromAJM1941}.

 If $p$ is an irreducible element of the monoid $\CM$ then so are
the elements $up$ and $pu$ for all
 units $u\in \CM^*$.
\begin{itemize} \item  {\em Each element $a$ of $\CM$ with $|a|>1$
is a product of irreducible elements.} \end{itemize}
 To prove this
statement we use induction on $|a|$. By (A2) and (A3), each
element $a$ with $|a|=2$ is irreducible. Suppose that $|a|>2$ and
the result holds for all elements $a'$ of $\CM$ with $1<|a'|<|a|$.
Then either the element $a$ is irreducible or, otherwise, it  is a
product, say  $bc$,  of two non-units $b$ and $c$. Since $|a|=
|b|\, | c|$, $|b|>1$ and $|c|>1$ (see (A2) and (A3)), we have
$1<|b|<|a|$ and $1<|c|<|a|$. By induction, the elements $b$ and
$c$ are products of irreducible elements, then so is the element
$a$. $\Box$

\begin{corollary}\label{1Sep7}%\marginpar{1Sep7}
Let $\CM$ be a reduction monoid, $p$ and $q$ be irreducible
elements of $\CM$ such that $\CM^* p \neq \CM^*q$ and there exists
an element $a\in \CM p \cap \CM q$ with $ |a|>1$. Then the natural
numbers $|p|$ and $|q|$ are co-prime.
\end{corollary}

{\it Proof}. Suppose that the natural numbers $|p|$ and $|q|$ are
not co-prime, i.e. $k:= \gcd (|p|, |q|)>1$, we seek a
contradiction. Then $|p|= ki$, $|q|=kj$ for some co-prime natural
numbers $i$ and $j$. By (A5), $\CM p \cap \CM q = \CM c$ for some
element $c$ of $\CM$ with $|c|= \lcm (|p|, |q|)=ijk$. Then
$c=\alpha p = \beta q$ for some elements $\alpha $ and $\beta$ of
$\CM$ with $|\alpha | = j$ and $|\beta | = i$. By (A6), there
exist elements $p_1, q_1, d\in \CM$ such that $p= p_1d$, $q=q_1d$,
$|d|=k>1$, $|p_1|= i$, $|q_1|= j$.

If $i=j=1$ then $|\alpha |= |\beta | =1$, and so $\alpha , \beta
\in \CM^*$, by (A3). The equality $ \alpha p = \beta q$ implies
that $\CM^* p = \CM^* q$. This contradicts to the assumption of
the corollary.

Therefore, either $i>1$ or $j>1$ or both $i$ and $j$ are strictly
greater than 1. These mean that either the element $p$ is
reducible (since $p=p_1d$, $|p_1|=i>1$, $|d|>1$) or the element
$q$ is reducible (since $ q=q_1d$, $|q_1|=j>1$, $|d|>1$) or both
elements $p$ and $q$ are reducible. These contradictions prove the
corollary. $\Box $

$\noindent $

{\bf Proof of Theorem \ref{b31Aug7}}.

$\noindent $

 1. The first statement is an
easy corollary of the second (since in the case (a):
$|p_iu|=|p_i|$ and $|u^{-1}p_{i+1}|= |p_{i+1}|$, by (A2) and
(A3)).

2. For each element $b$ of the monoid $\CM$ with $|b|>1$, let
$\Dec (b)$ be the set of all decompositions of the element $b$
into irreducible components. Two such decompositions, say $X$ and
$Y$, are equivalent, $X\sim Y$, if one can be produced from the
other by finitely many transformations of the types (a) and (b).
Clearly, this is an equivalence relation on the set $\Dec (b)$.
Let $X,Y\in \Dec (b)$ and $X', Y'\in \Dec (b')$. If $X\sim  Y$
then $ XX'\sim YX'$ in $\Dec (bb')$ and $X'X\sim X'Y$ in $\Dec
(b'b)$. If $X\sim Y$ and $X'\sim Y'$ then $XX'\sim YY'$ in $\Dec
(bb')$.

To finish the proof of statement 2 we have to show that $p_1\cdots
p_r\sim q_1\cdots q_s$. To prove this fact we use induction on
$|a|$. Note that if the element $a$ is irreducible then $\Dec (a)=
\{ a\}$, and there is nothing to prove. The base of the induction,
$|a|=2$, is obvious since the element $a$ is irreducible, by (A2)
and (A3). Suppose that $|a|\geq 3$ and the result is true for all
elements $a'$ of $\CM$ with $1<|a'|<|a|$. We may assume that the
element $a$ is reducible, i.e. $r\geq 2$ and $s\geq 2$. The proof
consists of considering several possibilities.

Suppose that $\CM^* p_r= \CM^*q_s$, i.e. $p_r= uq_s$ for some
element $u\in \CM^*$. By (A4), we can delete the element $q_s$ in
the equality
$$ p_1\cdots p_{r-1}uq_s= q_1\cdots q_{s-1}q_s.$$
As a result, there are two decompositions of the element
$$a':=p_1\cdots p_{r-1}u= q_1\cdots q_{s-1}$$
into irreducible components with $1< |a'|= \frac{|a|}{|q_s|}<|a|$
(note that $p_{r-1}u$ is the irreducible element). By induction,
these two decompositions are equivalent in $\Dec (a')$. In
particular, $r=s$. Now,
$$ p_1\cdots p_r\sim p_1\cdots (p_{r-1}u)(u^{-1} p_r)\sim p_1\cdots (p_{r-1}u)\cdot q_s\sim q_1\cdots q_{r-1}\cdot
q_s,$$
 as required.

Suppose that $\CM^*p_r\neq \CM^*q_s$. Then, by Corollary
\ref{1Sep7}, the natural numbers $|p_r|$ and $|q_s|$ are co-prime
since $a=p_1\cdots p_r=q_1\cdots q_s\in \CM p_r\cap \CM q_s$ and the
elements $p_r$ and $q_s$ are irreducible. By (A6),
$$ \CM p_r\cap \CM q_s = \CM c$$
 for some element $c$ of the monoid  $\CM$ with $|c|= \lcm (
 |p_r|, |q_s|)= |p_r||q_s|$ since the numbers $|p_r|$ and $|q_s|$
 are co-prime. Since $a\in \CM c$ and $ c\in \CM p_r\cap \CM q_s$,
 there exist elements $d, \alpha , \beta \in \CM$ such that
 %\marginpar{dar1}
\begin{equation}\label{dar1}
a=dc, \;\; c= \alpha p_r= \beta q_s.
\end{equation}
We can write the equality $dc= a$ in two different ways:
$$ d\alpha p_r= p_1\cdots p_{r-1}p_r  \;\; {\rm and } \;\; d\beta q_s=
q_1\cdots q_{s-1}q_s.$$ By (A4), we can delete the element $p_r$
in the first equality and the element $q_s$ in the second:
%\marginpar{dar2}
\begin{equation}\label{dar2}
d\alpha = p_1\cdots p_{r-1} \;\; {\rm and } \;\; d\beta =
q_1\cdots q_{s-1}.
\end{equation}
Note that $1<|p_1|\leq |d\alpha |= \frac{|a|}{|p_r|}<|a|$ and $1<
|q_1|\leq |d\beta |= \frac{|a|}{|q_s|}<|a|$ since $r,s\geq 2$.
Then induction yields the equivalence relations
$$ d\alpha \sim p_1\cdots p_{r-1}\;\; {\rm and}\;\; d\beta \sim
q_1\cdots q_{s-1}.$$ There are two options: either $|d|>1$ or
$|d|=1$.

If $|d|>1$ then $1< |p_r|\leq |c| = \frac{|a|}{|d|}<|a|$ (see
(\ref{dar1})), and so, by induction,  $\alpha p_r\sim \beta q_s$.
Now,
$$ p_1\cdots p_{r-1}p_r\sim d\alpha p_r\sim d\beta q_s\sim
q_1\cdots q_{s-1} q_s.$$

Finally, suppose that $|d|=1$. By (A3), the element $d$ is a unit of
the monoid $\CM$ since $|d|=1$. Then $\CM c = \CM da = \CM a$ (since
$c = da$). Without loss of generality we may assume that $c=a$ and
$d=1$. Then the equations (\ref{dar2}) mean that %\marginpar{dar3}
\begin{equation}\label{dar3}
\alpha = p_1\cdots p_{r-1}\;\; {\rm and }\;\; \beta = q_1\cdots
q_{s-1}.
\end{equation}
Recall that we have the equality $|c|= |p_r||q_s|$. In combination
with (\ref{dar1}), i.e. $a= c = \alpha p_r = \beta q_s$, it yields
the equalities
 $$ |\alpha | = |q_s| \;\; {\rm and } \;\; |\beta | = |p_r|.$$
In particular, the numbers $|\alpha |$ and $|\beta |$ are
co-prime. Recall that $r\geq 2$ and $s\geq 2$. Now, the case
$r=s=2$ is trivially true, $p_1p_2\sim q_1q_2$, since $a=p_1p_2=
q_1q_2$ and the numbers $|p_1|= |q_2|$ and $|p_2|=|q_1|$ are
co-prime. This is a transformation of the type (b).

It remains to consider the case $(r,s)\neq (2,2)$.  In a view of
symmetry, we may assume that $r\geq 3$ and $s\geq 2$. We prove that
this case is not possible, i.e. we seek a contradiction. In order to
get a contradiction, the axiom (A6) will be applied to the equality
%\marginpar{dar4}
\begin{equation}\label{dar4}
p_1\cdot (p_2\cdots p_r)= \beta \cdot q_s.
\end{equation}
First, note that the numbers
$$ i:= |p_1|= \frac{|p_1\cdots p_{r-1}|}{|p_2\cdots p_{r-1}|}=\frac{|\alpha |}{|p_2\cdots
p_{r-1}|}=\frac{|q_s|}{|p_2\cdots p_{r-1}|}\;\; {\rm and}\;\; j:=
|\beta |= |p_r|$$ are co-prime since the  numbers $|q_s|$ and
$|p_r|$ are co-prime; $i>1$ and $j>1$. Clearly, $k:= |p_2\cdots
p_{r-1}|>1$ since $r\geq 3$; $|p_2\cdots p_r|=kj$ and $|q_s|= ki$.
Applying the axiom (A6) to the equality (\ref{dar4}), we obtain the
equalities
$$ p_2\cdots p_r=AC \;\; {\rm and}\;\; q_s= BC$$
for some elements $A$, $B$ and $C$ of the monoid $\CM$ with $|C|=
k>1$. Then $|B|= \frac{|q_s|}{|C|}=\frac{ki}{k}=i>1$, and so the
elements $B$ and $C$ are not units. Therefore, the element $q_s=
BC$ is reducible, a contradiction. The proof of Theorem
\ref{b31Aug7} is complete. $\Box$

$\noindent $

{\bf Proof of Theorem \ref{2Nov7}.}

$\noindent $

In the proof of Theorem \ref{2Nov7}, we use the Theorem of ${\rm
L}\ddot{{\rm u}}{\rm roth}$  and the fact that $\OO$ is a submonoid
of the reduction monoid $(K[x], \circ )$. The axioms (A1)--(A4) are
obvious for the monoid $\OO$.

Let us prove that the axiom (A5) holds for $\OO$. Let $a$ and $b$ be
elements of the monoid $\OO$ such that $\deg (a)>1$, $\deg (b)>1$,
and there exists an element $x'\in \OO \circ a \cap \OO \circ b$
with $\deg (x')\geq 1$. Note that $x'\in \OO$. Then $x'\in K[x]\circ
a \cap K[x]\circ b$, and so $K[x]\circ a \cap K[x]\circ b=K[x]\circ
c$ for some element $c$ of $K[x]$, by the axiom (A5) for the
reduction monoid $K[x]$. Moreover, $\deg (c) = \lcm (\deg (a) , \deg
(b))$.

It suffices to show that $c+\nu \in \OO$ for some element $\nu \in
K$. For, we introduce the $K$-algebra automorphism $\o$  of the
polynomial algebra $K[x]$
given by the rule $x\mapsto -x$. Then %\marginpar{Kx2CO}
\begin{equation}\label{Kx2CO}
K[x]= K[x^2]\oplus K[x^2]x=K[x^2]\oplus  \OO ,
\end{equation}
where $K[x^2]$ is the fixed ring for the automorphism $\o$, and
$\OO$ is the eigen-space for $\o$ that corresponds to the eigenvalue
$-1$, i.e. $\OO = \ker (\o +1)$.  Note that the equality $K[x]\circ
a \cap K[x]\circ b=K[x]\circ c$ simply means that
$$K[a] \cap K[b]= K[c],$$
and so the element $c$ is uniquely defined up to an affine
transformation. By (\ref{Kx2CO}), the element $c$ is a unique sum
$c_0+c_1x$ for some elements $c_0, c_1\in K[x^2]$.  Note that
$c_1\neq 0$ since, otherwise, $c= c_0\in K[x^2]$, and then
$$ x'\in \OO\circ a \cap \OO \circ b \subseteq K[x] \circ a \cap
K[x] \circ b = K[c]\subseteq K[x^2].$$ Now, $x'\in \OO \cap K[x^2]=
0$,  a contradiction  (recall that $\deg (x') \geq 1$, by the
assumption). This contradiction proves the claim that $c_1\neq 0$.
Note that
$$ \o (K[c])= \o (K[a]\cap K[b]) = \o (K[a]) \cap \o (K[b]) =
K[-a]\cap K[-b]= K[a]\cap K[b] = K[c].$$ This means that $\o (c) =
\l c +\mu$ for some scalars $\l \neq 0$ and $\mu$ of $K$. In
combination with the equality $\o (c) = c_0-c_1x$ and the fact that
$c_1\neq 0$, it gives that $\l = -1$,  i.e. $\o (c) = -c+\mu$. Then
changing $c$ to $c-\frac{\mu}{2}$ we may assume that $\mu =0$, i.e.
$\o (c) = -c$.  This means that $c\in \OO$, as required. This proves
that the axiom (A5) holds for the monoid $\OO$.

To finish the proof of Theorem \ref{2Nov7}, it remains to establish
the axiom (A6) for the monoid $\OO$.

Suppose that elements $a$, $b$, $\alpha$ and $\beta$ of the monoid
$\OO$ satisfy the following conditions: $\alpha \circ  a = \beta
\circ b$ with $\deg (\alpha )=i$, $\deg (a)= jk$, $\deg (\beta ) =
j$, $\deg (b)=ik$, $ ijk\geq 1$, and the natural numbers $i$ and $j$
are co-prime. We have to show that $a= a_1\circ d$ and $b= b_1\circ
d$ for some elements $a_1$, $b_1$ and $d$ of the monoid $\OO$ such
that $\deg (d) =k$. In the proof  of the axiom (A5) for the monoid
$\OO$, we found the element $c\in \OO$ such that
$$ K[c]= K[a]\cap K[b], \;\; \deg (c) = \lcm (\deg (a) , \deg (b)) =
ijk.$$ Then, it is easy to show that %\marginpar{Kc=KaKb}
\begin{equation}\label{Kc=KaKb}
K(c) = K(a) \cap K(b).
\end{equation}
Indeed, by the Theorem of ${\rm L}\ddot{{\rm u}}{\rm roth}$, $ K(a)
\cap K(b)= K(c^*)$ for some element $c^*\in K(x)$ which can be
chosen from the polynomial algebra $K[x]$, by Lemma 3.1,
\cite{EngstromAJM1941}. Then
$$ K[c^*] = K[x]\cap K(c^*)= (K[x] \cap K(a))\cap (K[x]\cap K(b)) =
K[a]\cap K[b]= K[c],$$ and so the equality (\ref{Kc=KaKb}) follows.

For a field extension $\D \subseteq \G$, let $[\G : \D ] := \dim_\D
(\G )$. Consider the fields $K(c)\subseteq K(a) \subseteq K(x)$.
Then
\begin{eqnarray*}
 ijk&=& \deg (c) = [K(x):K(c)]= [ K(x):K(a)]\cdot [ K(a):K(c)] \\
 &=& \deg (a) \cdot [ K(a):K(c)] =jk\cdot [ K(a):K(c)],
\end{eqnarray*}
hence $[K(a):K(c)]=i$. By symmetry, $[ K(b):K(c)]=j$. By the Theorem
of ${\rm L}\ddot{{\rm u}}{\rm roth}$, the composite field $K(a) K(b)
= K(a,b)\subseteq K(x)$ is equal to $K(d)$ for some rational
function $d\in K(x)$ which can be chosen to be a polynomial of
$K[x]$ since $a,b\in K[x]$. Let us show that %\marginpar{Kd=ij}
\begin{equation}\label{Kd=ij}
[K(d):K(c)]=ij.
\end{equation}
Clearly,
\begin{eqnarray*}
 [K(d):K(c)]&=&  [K(a,b):K(c)]=[K(a)(b):K(a)] [ K(a):K(c)]\\
 &\leq & [K(c)(b):K(c)] [ K(a):K(c)] \\
 &=& [ K(b):K(c)]  [ K(a):K(c)]=ji.
\end{eqnarray*}
To prove the reverse inequality note that
$$[K(d):K(c)]=  [K(d):K(a)] [ K(a):K(c)]=[K(d):K(a)]\cdot i, $$
$$ [K(d):K(c)]=  [K(d):K(b)] [ K(b):K(c)]=[K(d):K(b)]\cdot j,$$
and so $[K(d):K(c)]\geq \lcm (i,j) = ij$ since the numbers $i$ and
$j$ are co-prime. This proves the equality (\ref{Kd=ij}). Now,
$$ \deg (d) = \frac{[K(x):K(c)]}{[K(d):K(c)]}=\frac{ijk}{ij}=k.$$
Note that
$$
 K(\o (d))= \o (K(d))= \o (K(a,b)) = K(\o (a) , \o (b))= K(-a, -b) = K(a,b) = K(d).
 $$
This means that $\o (d) = \l d +\mu$ for some scalars $\l \neq 0$
and $\mu$ of $K$ since $d\in K[x]$ and $\o (K[x]) = K[x]$. By
(\ref{Kx2CO}), the polynomial $d$ is a unique sum $d_0+d_1x$ for
some polynomials $d_0, d_1\in K[x^2]$. We must have $d_1\neq 0$
since, otherwise, $d=d_0\in K[x^2]$. Since $a=a_1\circ d$ for some
polynomial $a_1\in K[x]$, we would have $a\in a_0\circ
K[x^2]\subseteq K[x^2]$, and so $a\in \OO\cap K[x^2]= 0$, a
contradiction (since $a \neq 0$). Therefore, $d_1\neq 0$. Then the
equalities
$$ d_0-d_1x= \o (d) = \l d +\mu = \l d_0+\mu +\l d_1x$$
yield $\l =-1$, and so $\o ( d) =   -d+\mu$. Then changing $d$ for
$d-\frac{\mu}{2}$ we may assume that $\mu =0$, that is $\o (d) =
-d$,  i.e. $d\in \OO$. We claim that the polynomial $a_1\in K[x]$ in
the equality $a= a_1\circ d$ above belongs to $\OO$. To prove this
we write the polynomial $a_1$ as a unique sum $u+vx$ for some
polynomials $u,v\in K[x^2]$. Note that $u\circ d, v\circ d\in
K[x^2]$ and $(v\circ d) \cdot d\in \OO$. The inclusion
$$a=a_1\circ d = u\circ d +(v\circ d) \cdot d \in \OO$$
yields $u\circ d=0$, i.e. $u=0$. This proves that $a_1=vx \in \OO$.
By symmetry, we have $b=b_1\circ d$ for some element $b_1\in \OO$.
This means that the axiom (A6) holds for the monoid $\OO$. The proof
of Theorem \ref{2Nov7} is complete. $\Box$

$\noindent $

{\it Definition}. A {\em Ritt transformation}  of the
decomposition (\ref{ap1p1}) is either one of the decompositions
(a), (b) or (c) with $\l_2=1$ and $\gcd (\deg (p_i), \deg
(p_{i+1}))=1$ (in all three cases) and with the numbers $m$ and
$n$ being {\em odd} prime numbers in the case (a) (see
(\ref{ap1p2})) or a decomposition of the type
\begin{eqnarray*}
{\rm (d)} & p_1\circ \cdots \circ (p_i\circ u)\circ (u^{-1}
\circ p_{i+1})\circ \cdots \circ p_r \\
\end{eqnarray*}
 for some polynomial $u\in
K[x]^*$.

$\noindent $

In his paper, J. F. Ritt wrote (page 52, the last line): {\em
``Case (a) with $m=2$ can be reduced to Case (b) by linear
transformation.''} In more detail, for each natural number $k\geq
1$,
\begin{eqnarray*}
 T_2&=& -1+2x^2= (-1+2x)\circ x^2=\alpha\circ x^2, \;\; \alpha :=-1+2x, \\
 T_{2k+1}&=& \sum_{i=0}^k{2k+1\choose 2i} x^{2k+1-2i}(1-x^2)^i=
 xt_{2k+1}(x^2), \\
 t_{2k+1}(x)&:=& \sum_{i=0}^k{2k+1\choose 2i} x^{k-i}(1-x)^i.
\end{eqnarray*}

Let $n=2k+1$. Then
$$ T_2\circ T_n=\alpha \circ x^2\circ [ xt_n(x^2)]=\alpha \circ [
xt_n^2]\circ x^2= \alpha \circ [xt_n^2]\circ \alpha^{-1}\circ
\alpha \circ x^2= \alpha \circ [ xt_n^2]\circ \alpha^{-1}\circ
T_2,
$$ and the remark of J. T. Ritt is obvious. Note that $T_n\circ
T_2= T_2\circ T_n = \alpha \circ [xt_n^2]\circ \alpha^{-1} \circ
T_2$, and so (by (A4))
$$ T_n= \alpha \circ [ xt_n^2]\circ \alpha^{-1}.$$
Now, it is obvious that also the case (a) with $n=2$ can be reduced
to the case (c) by linear transformation. This is the reason why in
the definition of Ritt transformation $m$ and $n$ are odd primes (in
the case (a)).

$\noindent $

All trigonometric polynomials $T_l= xt_l(x^2)$ {\em do not} belong
to the algebra $A$  where $l$ runs through all {\em odd} prime
numbers (since $T_l'(0)= l\neq 0$). But $T_2\in A$.

The next corollary  follows from Theorem \ref{b31Aug7} and the
second theorem of Ritt(-Levi), it is implicit in the papers
\cite{RittTrAMS1922} and \cite{LeviAJM1942}.

\begin{corollary}\label{c31Aug7}%\marginpar{c31Aug7}
If $a\in K[x]$ has two decompositions into irreducible polynomials
then one can be obtained from the other by Ritt transformations.
\end{corollary}

$\noindent $

{\bf Proof of Corollary \ref{a2Nov7}.(3).}

$\noindent $

The idea of the proof of Corollary \ref{a2Nov7}.(3) is to use the
second theorem of Ritt-Levi in combination with Lemma \ref{3Nov7},
Theorem \ref{8Nov7} and Lemma \ref{c8Nov7}. We first prove all these
preliminary results that are interesting on their own.

\begin{lemma}\label{3Nov7}%\marginpar{3Nov7}
Let $K$ be a field of characteristic zero, $a$ and $b$ be non-scalar
polynomials of $K[x]$ such that $a\circ b\in \OO$. If one of the
polynomials $a$ or $b$ belongs to   the set $\OO$ then so  does the
other.
\end{lemma}

{\it Proof}. {\em Case (i)}: $a\in \OO$. The polynomial $a$ is a
non-scalar polynomial, and so
$$ a= \sum_{n=0}^N\l_nx^{2n+1}, \;\; \l_n\in  K, \;\; \l_N\neq 0. $$
Due to the decomposition $K[x]= K[x^2]\oplus K[x^2]x$, each
polynomial $p$ of $K[x]$ is a unique sum $p = p^{ev}+p^{od}$ of an
even $p^{ev}\in K[x^2]$ and odd $p^{od}\in K[x^2]x$ polynomials.
Then  $b=b_0+b_1$ where $b_0:= b^{ev}$ and $b_1:= b^{od}$. We have
to show that $b_0=0$. Suppose that $b_0\neq 0$, we seek a
contradiction. Clearly, $b_1\neq 0$   since otherwise we would have
 the inclusion $c\in  K[x^2]x\cap K[x^2]=0$, a contradiction. Let us
 consider the even part  of the polynomial $c$,
 $$ c^{ev} =  (a\circ b)^{ev} =   (\sum_{n=0}^N \l_n
 (b_0+b_1)^{2n+1})^{ev}= \sum_{n=0}^N\l_n\sum_{m=0}^n{2n+1\choose
 2m+1}b_0^{2m+1}   b_1^{2(n-m)}.$$
The degrees of the nonzero polynomials $b_0$ and $b_1$ are even and
odd numbers respectively. Therefore,  either $\deg (b_0)>\deg (b_1)$
or, otherwise, $\deg (b_0)<\deg (b_1)$. The leading coefficient of
the polynomial $c^{ev}$ is equal  to
$$
\begin{cases}
\l_Nb_0^{2N+1}& \text{if $\deg (b_0)>\deg (b_1)$},\\
\l_N{2N+1\choose 1} b_0b_1^{2N}& \text{if $\deg (b_0)<\deg (b_1)$}.\\
\end{cases}
$$
The first case is obvious; the second case follows from  the
inequalities: for  all natural numbers $m$ and $n$ such that $0\leq
m \leq n$,
$$ \deg (b_0^{2m-1}b_1^{2(n-m+1)})-\deg
(b_0^{2m+1}b_1^{2(n-m)})=2(\deg (b_1)-\deg (b_0))>0.$$ Since in both
cases the leading term of the polynomial $c^{ev}$ is non-zero, we
have $c^{ev}\neq 0$. This contradicts to the assumption that $c\in
K[x^2]x$, i.e. $c^{ev}=0$. The contradiction  finishes the proof of
the case (i).

{\em Case (ii)}: $b\in \OO$. Then  $\o (b)=-b$. Similarly, $\o (c)
=-c$ since $c\in K[x^2]x$. The polynomial $a$ is a unique sum
$a^{ev}+a^{od}$ of even and odd polynomials. Comparing both ends of
the following series of equalities
\begin{eqnarray*}
-(a^{ev}\circ b+a^{od}\circ b)&=& -c=\o (c) = \o (a\circ b)=a\circ
\o
(b) = a\circ (-b)\\
&=& a^{ev}\circ b - a^{od}\circ b
\end{eqnarray*}
we conclude that $a^{ev}\circ b=0$, hence $a^{ev}=0$ since $b$ is a
non-scalar polynomial, and so $a=a^{od}\in \OO$, as required. The
proof of Lemma \ref{3Nov7} is complete. $\Box $

$\noindent $

Let $f=f_0+f_1\in K[x]$ where $f_0:=f^{ev}$ and $f_1:=f^{od}$. Let
$f^{(k)}:=\frac{d^kf}{dx^k}$ and $f^{(k)}(g):=\frac{d^kf}{dx^k}\circ
g$. Then $f^{(2n)}=f_0^{(2n)}+f_1^{(2n)}$ and
$f^{(2n+1)}=f_1^{(2n+1)}+f_0^{(2n+1)}$ where $f_0^{(2n)},
f_1^{(2n+1)}\in K[x^2]$ and $f_1^{(2n)}, f_0^{(2n+1)}\in \OO$.

\begin{lemma}\label{a8Nov7}%\marginpar{a8Nov7}
Let $f=f^{ev}+f^{od}\in K[x]$ and $\mu \in K^*$. Then $(x+\mu )
\circ f\in \OO$ iff $f^{ev}=-\mu$.
\end{lemma}

{\it Proof}. $(x+\mu ) \circ f=\mu +f=\mu +f^{ev}+f^{od}\in \OO$ iff
$f^{ev}=-\mu$. $\Box $

$\noindent $

\begin{lemma}\label{b8Nov7}%\marginpar{b8Nov7}
Let $a\in \OO$ and $f=f_0+f_1\in K[x]$ where $f_0:= f^{ev}$ and
$f_1:=f^{od}$. Then $(a\circ f)^{ev}= \sum_{k\geq 0}
a^{(2k+1)}(f_1)\cdot \frac{f_0^{2k+1}}{(2k+1)!}$ and $(a\circ
f)^{od}= \sum_{k\geq 0} a^{(2k)}(f_1)\cdot \frac{f_0^{2k}}{(2k)!}$.
\end{lemma}

{\it Proof}. The result is an easy consequence of the Taylor's
formula, $$a\circ f= a(f_1+f_0)= \sum_{i\geq 0} a^{(i)}(f_1)\cdot
\frac{f_0^i}{i!}, $$ and the following two facts:
$a^{(2k+1)}(f_1)\in K[x^2]$ and $a^{(2k)}(f_1)\in \OO$. $\Box $

$\noindent $

\begin{theorem}\label{8Nov7}%\marginpar{8Nov7}
Suppose that $a\in \OO$ with $\deg (a)>1$, $\mu \in K^*$, and $f\in
K[x]\backslash K$. Then $(x+\mu ) \circ a \circ f\not\in \OO$.
\end{theorem}

{\it Proof}. Suppose that $(x+\mu ) \circ a \circ f\in \OO$, we seek
a contradiction. Then
\begin{eqnarray*}
 -\mu &=& (a\circ f)^{ev} \;\;\;\;\;\;\;\;\;\; ({\rm by \; Lemma}\;  \ref{a8Nov7}) \\
 &=&\sum_{k\geq 0}a^{(2k+1)}(f_1)\cdot \frac{f_0^{2k+1}}{(2k+1)!} \;\;\;\; ({\rm by \; Lemma}\;  \ref{b8Nov7}) \\
 &=& f_0\cdot \sum_{k\geq 0}a^{(2k+1)}(f_1)\cdot
 \frac{f_0^{2k}}{(2k+1)!}.
\end{eqnarray*}
Comparing the degrees of both ends of the series of equalities
above, we conclude that $f_0\in K^*$ since $\mu \neq 0$. Let $\der
:= \frac{d}{dx}$. Then $-\mu = \D \der (f_1)$ where the linear map
$$\D := \sum_{k\geq 0} \frac{f_0^{2k+1}}{(2k+1)!}\der^{2k}:K[x]\ra
K[x]$$ is equal to $f_0 (1-\gn )$ where $\gn := -\sum_{k\geq 1}
\frac{f_0^{2k}}{(2k+1)!}\der^{2k}$ is a locally nilpotent map, that
is $K[x] = \cup_{i\geq 1} \ker (\D^i)$. The map $\D$ is invertible
and $\D^{-1} =f_0^{-1} (1+\gn + \gn^2+\cdots )$. Then $\der (f_1) =
-\D^{-1} (\mu ) = -f_0^{-1}$, and so $\deg (f_1)\leq 1$, that is
$f_1= \g x$ for some $\g \in K^*$ since $f\not\in K$. We claim that
$f_0\neq 0$ since otherwise we would have the inclusion $ (x+\mu )
\circ a \circ f_1= a\circ f_1+\mu \in \OO$, which would have implied
that $\mu =0$ (since $a\circ f_1\in \OO$), a contradiction.
Changing, if necessary, the element $a$ to $a\circ f_1 = a\circ [\g
x]\in \OO$, we may assume that $\g =1$. Then $\OO \ni (x+\mu ) \circ
a \circ (f_0+x)$ iff
$$ -\mu = (a\circ (f_0+x))^{ev} = \sum_{k\geq 0}
a^{(2k+1)}\frac{f_0^{2k+1}}{(2k+1)!} \;\;\;\; ({\rm see \;
above}).$$ This implies that $\deg (a) \leq 1$, a contradiction
(since $\deg (a) >1$, by the assumption). This contradiction
finishes the proof of the theorem. $\Box $

$\noindent $

The next corollary follows at once from Theorem \ref{8Nov7}.
\begin{corollary}\label{x8Nov7}%\marginpar{x8Nov7}
Suppose that $a\in \OO$ with $\deg (a) \geq 1$, $\mu \in K^*$, and
$f\in K[x]\backslash K$. If $(x+\mu ) \circ a \circ f\in \OO$ then
$\deg (a)=1$.
\end{corollary}

$\noindent $

{\it Example}. $(x+\mu ) \circ [ \l x]\circ (x-\l^{-1} \mu )\in \OO$
for all $\l , \mu \in K^*$.

$\noindent $

\begin{lemma}\label{c8Nov7}%\marginpar{c8Nov7}
Let $f\in K[x]$ with $\deg (f) \geq 1$ and $u \in K[x]^*$. Then
$f\circ x^2\circ u \not\in \OO$.
\end{lemma}

{\it Proof}. Let $u=\l x+\mu$ for some $\l \in K^*$, and
$f=\sum_{i=0}^n\l_ix^i$ where $n:=\deg (f)$, and so $\l_n\neq 0$.
Then $f\circ x^2\circ (\l x+\mu )=\sum_{i=0}^n (\l x+\mu )^{2i} =
\l_n \l^{2n}x^{2n}+{\rm smaller\;  terms}$, and so $f\circ x^2\circ
u \not\in \OO$.  $\Box $

$\noindent $

{\bf The proof of Corollary \ref{a2Nov7}.(3) continued}.  Let us
continue with the proof of Corollary \ref{a2Nov7}.(3). Recall that
$\OO^*= \{ \l x\, | \, \l \in K^*\}$. We have to show that if there
is an equality $p\circ q= p^*\circ q^*$ where $p$, $q$, $p^*$ and
$q^*$ are irreducible elements of the monoid $\OO$ then  modulo {\em
basic} transformations of the pairs $P:=(p,q)$ and $P^*:= (p^*,
q^*)$:
$$(p,q)\mapsto (u\circ p \circ v , v^{-1}\circ  q\circ  w), \;\;
(p^*,q^*)\mapsto (u\circ p^* \circ \widetilde{v} ,
\widetilde{v}^{-1}\circ q\circ w), $$ where $u$, $v$,
$\widetilde{v}$, $w\in \OO^*$, we have  either the equality $P=P^*$
 or, otherwise, $P$ and $P^*$ as in Corollary \ref{a2Nov7}.(3).

If $(p^*, q^*)=(p\circ v, v^{-1}\circ q)$ for some element $v\in
K[x]^*$ then, by Lemma \ref{3Nov7}, $v\in \OO^*$, and there is
nothing to prove, the result is obvious. So, suppose that $(p^*,
q^*)\neq (p\circ v, v^{-1}\circ q)$ for all
 element $v\in K[x]^*$. Then  by the second theorem of Ritt-Levi the
 pair $P^*$ can be   obtained from the pair $P$ by finitely many
 Ritt transformations
 $$P=P_1\sim_RP_2\sim_R\cdots \sim_R P_s=P^*,$$
 and necessarily some of the Ritt transformations are of the types (a),
 (b) or (c). It might happen that the elements $p$ and $q$ are reducible
 in the monoid $K[x]$ (but the essence of the proof is to show that
 they are, in fact, irreducible in $K[x]$).

Each Ritt transformation $P_i:= (p_i, q_i) \sim_R P_{i+1}:=(p_{i+1},
q_{i+1})$ may transform either the irreducible factors (in $(K[x],
\circ )$) of $p_i$ or of $q_i$ or simultaneously the last
irreducible factor, say $l_i$, of $p_i$ and the first irreducible
factor, say $f_i$, of $q_i$. The first  two types of Ritt
transformations do not change the elements $p_i$ and $q_i$. So,
there exists an index $i$ such that the Ritt transformation
$P_i\sim_R P_{i+1}$ is of the third type and, necessarily, of one of
the types (a), (b) or (c) as in the definition of Ritt
transformations since, for  given $u\in K[x]^*$ and $a\in \OO^*$,
the inclusion $u\circ a\in \OO^*$ implies $u\in \OO^*$ (Lemma
\ref{3Nov7}). Let $i$ be the least such an index. For each $j$, let
$Q_j:=(l_j, f_j)$. Then $p_j=\alpha_j\circ l_j$ and $q_j=f_j\circ
\beta_j$ for some polynomials $\alpha_j, \beta_j\in K[x]$. There are
the following three options for the pairs $Q_i=(l_i, f_i)$ and
$Q_{i+1}=(l_{i+1}, f_{i+1})$ (where $u,v, w, \widetilde{w}\in
K[x]^*$):

 (a) $Q_i=(u\circ T_n\circ w, w^{-1}\circ T_m\circ    v)$ and
$Q_{i+1}=(u\circ T_m\circ \widetilde{w}, \widetilde{w}^{-1}\circ
T_n\circ v)$ where $n$ and $m$ are odd primes,

 (b) $Q_i= (u\circ [
x^t\beta^s]\circ w, w^{-1}\circ x^s\circ v)$ and $Q_{i+1}=(u\circ
x^s\circ \widetilde{w}, \widetilde{w}^{-1}\circ  [ x^t\beta
(x^s)]\circ v)$,

(c)  $Q_i=(u\circ x^s\circ w,  w^{-1}\circ   [   x^t\beta
(x^s)]\circ v)$ and $Q_i= (u\circ [ x^t\beta^s]\circ \widetilde{w},
\widetilde{w}^{-1}\circ x^s\circ v)$,

where $s$ is a  prime number, $t\geq 0$,  and $\beta \in K[x]$ with
$\beta (0)\neq 0$.  In    the  cases (b) and (c), $s$ is an {\em
odd} prime number since, otherwise, by Lemma \ref{c8Nov7},  the
polynomials $p_{i+1}\not\in \OO$ (the case (b)) and $p_i\not\in \OO$
(the case (c)), which are contradictions.
% Since the polynomials
% $T_n$, $T_m$ and $x^s$ are elements of the monoid $\OO$, and the
% polynomials $p_j=\alpha_j\circ l_j$  and $q_j= f_j\circ \beta_j$ are
 % irreducible elements of the monoid $\OO$ for $j=i, i+1$, then, by Lemma
% \ref{3Nov7}, the polynomials $p_i$ and $q_i$ are {\em irreducible}
% polynomials of the monoid $K[x]$ in all three cases (a)--(c), i.e.
% $Q_i=(p_i, q_i)$. It follows directly from the fact that $\beta
% (0)\neq 0$ that the polynomial $x^t\beta (x^s)$ (see the cases (b)
% and (c)) belongs to the monoid $\OO$ iff $t$ is an odd number and
% $\beta (x^s) = \alpha (x^{2s})$ for some polynomial $\alpha (x) \in
% K[x]\backslash K$ with $\alpha (0)\neq 0$.

Let us consider the case (a). Note that $T_m, T_n\in \OO$. Applying
Theorem \ref{8Nov7} to the inclusion $w^{-1} \circ T_m\circ (v\circ
\beta_i)= q_i\in \OO$, we see that $w^{-1}\in \OO^*$. Then we have
the inclusion $T_m\circ (v\circ \beta_i)\in \OO$ which yields the
inclusion $v\circ \beta_i\in \OO$, by Lemma \ref{3Nov7} (since
$T_m\in \OO$). Since $q_i$ is an irreducible element of the monoid
$\OO$, we must have $v\circ \beta_i\in \OO^*$.

Since $w\in\OO^*$ and $(\alpha_i\circ u \circ T_n) \circ w =p_i\in
\OO$, we have the inclusion $\alpha_i\circ u \circ T_n\in \OO$,
hence $\alpha_i\circ u \in \OO$ (by Lemma \ref{3Nov7} since $T_n\in
\OO$). Moreover, $\alpha_i\circ u \in \OO^*$ since $p_i$ is an
irreducible element of the monoid $\OO$. As a result, we have the
case (a) of Corollary \ref{a2Nov7}.(3).

Let us consider the case (b). Since $x^s\in \OO$ and $w^{-1} \circ
x^s\circ (v\circ \beta_i) = q_i\in \OO$, we have $w^{-1} \in \OO^*$
(by Theorem \ref{8Nov7}). Then $x^s\circ (v\circ \beta_i) \in \OO$,
hence $v\circ \beta_i\in \OO$, by Lemma \ref{3Nov7}. The element
$q_i$ is an irreducible element of the monoid $\OO$, and so $v\circ
\beta_i\in \OO^*$. By replacing the element $v$ with $v\circ
\beta_i$, we may assume that $\beta_i=1$ and $v\in \OO^*$. Now, it
follows from the inclusion
$$ \OO \ni q_{i+1} = \widetilde{w}^{-1} \circ [x^t\beta(x^s)]\circ v
\circ \beta_i=\widetilde{w}^{-1} \circ [x^t\beta(x^s)]\circ v$$ that
$\widetilde{w}^{-1} \circ [x^t\beta(x^s)]\in \OO$.

If $t\neq 0$ then $\widetilde{w}^{-1}\in \OO^*$, and so
$x^t\beta(x^s)\in \OO$, hence $t$ is odd (since $\beta (0)\neq 0$),
and
 $\beta =\alpha (x^2)$ for some polynomial $\alpha (x) \in K[x]$.
Since $[x^t\beta(x^s)]\circ w\in \OO$ and $(\alpha_i\circ u)\circ
[x^t\beta(x^s)]\circ w=p_i\in \OO$, we have $\alpha_i\circ u \in
\OO$, by Lemma \ref{3Nov7}. Therefore, $\alpha_i\circ u \in \OO^*$
since $p_i$ is an irreducible  element of the monoid $\OO$ and
$x^t\beta (x^s) \not\in \OO^*$. This means that we have the case (b)
of Corollary \ref{a2Nov7}.(3) (if $t\neq 0$).

To finish with the case $b$ it suffices to show that the remaining
subcase when $t=0$ is impossible. Suppose that $t=0$, we seek a
contradiction. Then the inclusion $\widetilde{w}^{-1} \circ
\beta(x^s)\in \OO$ yields $\beta = \widetilde{w}\circ
x^T\alpha_1(x^2)$ for some odd natural number $T$ and a polynomial
$\alpha_1 (x) \in K[x]$  with $\alpha_1(0)\neq 0$. Note that $w\in
\OO^*$ and
$$ \OO \ni p_i= \alpha_i \circ u \circ \beta^s \circ w = \alpha_i
\circ u \circ x^s\circ \widetilde{w}\circ [ x^T\alpha_1(x^2)]\circ
 w.$$ Since $x^T\alpha_1(x^2)\circ w\in \OO$ and the element $p_i\in \OO$ is
irreducible, we must have $\alpha_i \circ u \circ x^s\circ
\widetilde{w}\in \OO^*$, by Lemma \ref{3Nov7}, hence $s=1$, a
contradiction ($s$ is a prime number).

The remaining case (c) follows from the case (b) by interchanging
the roles of the pairs (and repeating the proof of the case (b)).

 Therefore, the pairs
$P_i$ and $P_{i+1}$ are as in Corollary \ref{a2Nov7}.(3). By the
minimality of $i$, we have $p=p_1=\cdots = p_i$  and $q=q_1=\cdots =
q_i$, and so $P=P_i$. Now, the result is obvious. The proof of
Corollary \ref{a2Nov7}.(3)  is complete. $\Box$

$\noindent$

{\it Remark}. Let us explain the remark made in the Introduction
that the monoid $\OO$ has non-commutative origin. Let $\l$ be a
nonzero scalar. The algebra
$$ \L = \langle x,y \, | \, xy = \l yx \rangle$$ is called the {\em
quantum plane}. The algebra $\L$ is the skew polynomial algebra
$K[y][x; \s]$ where $\s$ is the $K$-algebra automorphism of the
polynomial algebra $K[y]$ which is given by the rule $\s (y) =\l y$.
The localization $\L':=S^{-1}\L$ of the algebra $\L$ at the Ore set
$S:= K[y]\backslash \{ 0 \}$ is the skew polynomial algebra $\L'=
K(y)[x;\s]$. Let $\l=-1$. The centre $Z'$ of the algebra $\L'$ is
the polynomial algebra $K(y^2)[x^2]$ with coefficients from the
field $K(y^2)$. Clearly,
$$\L'=K(y)[x^2]\oplus K(y)[x^2]x$$
where the algebra $K(y)[x^2]$ is the fixed ring of the inner
automorphism $\o_y:u\mapsto yuy^{-1}$ of $\L'$, and $K(y)[x^2]x=
\ker (\o_y+1)$. Then it follows that  the monoid $\CE$ of all the
$K$-algebra endomorphisms of $\L'$ elements of which fix the element
$y$ is equal to the set $\{ \tau_\alpha : x\mapsto \alpha x\, | \,
\alpha \in K(y) [ x^2]\}$. The endomorphism $\tau_\alpha$ is called
a {\em central} endomorphism if $\alpha \in Z'$. The submonoid $\CZ
:= \{ \tau_\alpha \, | \, \alpha \in Z'\}$ of all central
endomorphisms of $\L'$ is isomorphic to the monoid $\OO$ of odd
polynomials in $x$ where the base field is $K(y^2)$ rather than $K$.

$\noindent$

The set $\IrrK[x]$ of all the irreducible elements of the monoid
$(K[x], \circ )$ is the union of its three subsets,
%\marginpar{IPR}
\begin{equation}\label{IPR}
\IrrK[x] = \CP \cup \CQ \cup \CR
\end{equation}
where an irreducible polynomial $p$ is an element of the set $\CP$
iff $p\in K[x]^*\circ x^l\circ K[x]^*$ for some prime number $l$;
an irreducible polynomial $p$ belongs to $\CQ$ iff either
$$p\in
K[x]^*\circ [ x^sg(x^l)]\circ K[x]^*\;\; {\rm  or}\;\;  p\in
K[x]^*\circ [ x^sg^l]\circ K[x]^*$$
 for some prime number $l$, $s\geq 1$, $g(x) \in K[x]\backslash K$ with $g(0)\neq 0$; $\CR :=
\IrrK[x]\backslash \CP \cup \CQ$.

\begin{proposition}\label{a20Sep7}%\marginpar{a20Sep7}
\begin{enumerate}
\item The union (\ref{IPR}) is a disjoint union. \item The set
$\CP \cup \CQ$ contains precisely all the irreducible polynomials
of $K[x]$ that are involved in all the Ritt transformations.
\end{enumerate}
\end{proposition}

{\it Proof}. 1. By Lemma \ref{a19Sep7}, the union $\CP \cup \CQ$
is disjoint. Now, statement 1 is obvious.

2. For a prime number $l$, a polynomial $f$ of the form $g(x^l)=
g(x) \circ x^l$ (resp. $ g^l = x^l\circ g$) is irreducible iff $f\in
\CP$ (then, necessarily, $g$ is a unit). By   Lemma \ref{b19Sep7}
and the explicit formula for $T_l$ (see above), for each odd prime
number $l$,
$$ K[x]^* \circ T_l\circ K[x]^*\subseteq \CQ .$$
But $T_2\in \CP$. Now, statement 2 follows from the definitions of
Ritt transformations and of the sets $\CP$ and $\CQ$.  $\Box $

$\noindent $

\begin{lemma}\label{a19Sep7}%\marginpar{a19Sep7}
Let $f(x)$ be a non-scalar polynomial of $K[x]$ such that
$f(0)\neq 0$, $s$ and $p$ be natural numbers such that $s\geq 1$
and $p\geq 2$. Then the polynomials $x^sf(x^p)$ and $x^sf^p$ do
not belong to the set $\CN := \cup _{n\geq 2} K[x]^* \circ
x^n\circ K[x]^*$.
\end{lemma}

{\it Proof}. Suppose that $x^sf(x^p)\in \CN$, that is
$x^sf(x^p)=u\circ x^n\circ v$ for some elements $u$ and $v$ of the
set $K[x]^*$ and $n\geq 2$. We seek a contradiction. The derivative
$(u\circ x^n\circ v)'$ of the polynomial $u\circ x^n\circ v$ has a
single root with multiplicity $n-1\geq 1$. The same is true for the
derivative of the polynomial $x^sf(x^p)$ which is equal to
$$ (x^sf(x^p))'= x^{s-1}(sf(x^p)+px^pf'(x^p))=
x^{s-1}L(x^p)\neq 0$$ where $L(x) := sf(x) +pxf'(x)$. If $s\geq 2$
then zero must be a root of the polynomial $L(x^p)$, but $L(0) =
sf(0)\neq 0$, a contradiction. If $s=1$ then the polynomial
$L(x^p)$ must have a single root, say $\l $, which is not equal to
zero since $L(0)\neq 0$. Let $e$ be a $p$'th root of $1$ which is
not equal to $1$. Then $e\l$ is another root of $L(x^p)$ distinct
from $\l$, a contradiction. Therefore, $x^sf(x^p)\not\in \CN$.

Suppose that $x^sf^p(x)\in \CN$, that is $x^sf^p(x)=u\circ x^n\circ
v$ for some elements $u$ and $v$ of the set $K[x]^*$ and $n\geq 2$.
We seek a contradiction. By the same argument as in the previous
case, the derivative $(x^sf^p)'$ of the polynomial $x^sf^p$ must
have a single root with multiplicity $n-1\geq 1$. Clearly,
$$0\neq (x^sf^p)'= x^{s-1}\cdot f^{p-1}\cdot (sf+pxf').$$
Note that the polynomial $f^{p-1}$ has a nonzero root since
$f(0)\neq 0$. Hence, $s=1$ and the polynomials $f^{p-1}$ and
$f+pxf'$ have the same root, say $\l$,  but may be with different
multiplicities. The root $\l$ is a nonzero one since $f(0)\neq 0$.
Then $f= \mu (x-\l )^m$ for some $0\neq \mu \in K$ and $m\geq 1$,
and so
$$f+pxf'= \mu (x-\l )^{m-1} (x-\l + pmx).$$
Hence, $\l = \l (1+pm)^{-1}$, and so $1= 1+pm>1$, a contradiction.
Therefore, $x^sf^p(x)\not\in \CN$. $\Box $

$\noindent $

\begin{lemma}\label{b19Sep7}%\marginpar{b19Sep7}
Let $p$ be an odd natural number such that $p\geq 3$. Then the
trigonometric polynomial $T_p$ does not  belong to the set $\CN :=
\cup _{n\geq 2} K[x]^* \circ x^n\circ K[x]^*$.
\end{lemma}

{\it Proof}. The derivative $T_p'$ of the polynomial $T_p$ has at
least two distinct roots (Lemma \ref{RTkTl}) since $p\geq 3$, and
so the result. $\Box $

$\noindent $

The next result will be used in the proof of Theorem
\ref{A15Sep7}.
\begin{lemma}\label{RTkTl}%\marginpar{RTkTl}
Let $p$ be a natural number such that $p\geq 2$. Then
\begin{enumerate}
\item The derivative $T_p'$ of the trigonometric polynomial $T_p$
is a polynomial of degree $p-1$ which has $p-1$ distinct roots:
$\cos (\frac{\pi i }{p})$, $i=1, 2, \ldots , p-1$. \item If $k$
and $l$ are distinct prime numbers then the polynomials $T_k'$ and
$T_l'$ have no common roots.
\end{enumerate}
\end{lemma}

{\it Proof}. 1. By the very definition, the numbers $\cos
(\frac{\pi i }{p})$, $i=1, 2, \ldots , p-1$, are distinct. Note
that $\sin (\frac{\pi i }{p})\neq 0$ and  $\sin (p\cdot\frac{\pi i
}{p})=0$ for all $i=1, 2, \ldots , p-1$. Since $$T_p'(\cos (x))
\sin (x) = p \sin (px),$$ we have $T_p'(\cos (\frac{\pi i
}{p}))=0$ for all $i=1, 2, \ldots , p-1$. Now, statement 1 is
obvious since $\deg (T_p') = \deg (T_p) - 1 \leq p-1$.

2. Statement 2 follows from statement 1. $\Box $

$\noindent $

Let $a$ be a polynomial of $K[x]$ with $\deg (a)>1$ and $X= p_1\circ
\cdots \circ p_r\in \Dec (a)$ be a decomposition of the polynomial
$a$ into irreducible polynomials of $K[x]$. Let $n_\CP (X)$, $n_\CQ
(X)$ and $n_\CR (X)$ be  the numbers of irreducible factors $p_i$ of
the types $\CP$, $\CQ$ and $\CR$ respectively. For each prime number
$l$, let $n_{\CP , l}(X)$ be the number of irreducible factors $p_i$
such that $p_i\in K[x]^*\circ x^l\circ K[x]^*$.

\begin{theorem}\label{20Sep7}%\marginpar{20Sep7}
The numbers $n_\CP (X)$, $n_\CQ (X)$,  $n_\CR (X)$ and $n_{\CP ,
l}(X)$ do not depend on the decomposition $X$.
\end{theorem}

{\it Proof}. Recall that (\ref{IPR}) is a disjoint union, and the
set $\CP\cup \CQ$ contains precisely all the irreducible
polynomials that are involved in all the Ritt transformations
(Proposition \ref{a20Sep7}). Then it follows from the definition
of Ritt transformations that the numbers $n_\CP (X)$, $n_\CQ (X)$
and $n_{\CP , l}(X)$ do not depend on the decomposition $X$. Then
the number
$$ n_\CR = l(a) - n_\CP (X- n_\CQ (X)$$
 does not depend on the decomposition $X$ either. $\Box $

$\noindent $

{\it Definition}. The common value of all the numbers $n_\CP (X)$,
$X\in \Dec (a)$, is denoted by $n_\CP (a)$. Similarly, the numbers
$n_\CQ (a)$,  $n_\CR (a)$ and $n_{\CP , l}(a)$ are defined.

%%%%%%%%%%%%%%%%%% SECTION 3 %%%%%%%%%%%%%%%%%%%%%%%%

\section{Analogues of the two theorems of J. F. Ritt for
the  cusp}\label{CUSP}%\marginpar{CUSP}

In this section, Theorems \ref{15Sep7} and \ref{A15Sep7} are
proved. It is shown that, in general, the first theorem of J. F.
Ritt does not hold for the cusp, i.e., in general, the number of
irreducible polynomials in decomposition  of element of $A$ into
irreducible polynomials is not unique (Lemma \ref{b9Sep7}). For
each element $a$ of $A$, the set $\Max (a)$ is found (Lemma
\ref{22Sep7}).

In this section,   $K$ is an {\em algebraically closed}  field of
characteristic $0$ if it is not stated otherwise.

 The algebra $K[s,t]/(s^2-t^3)$ of regular
functions on the cusp $s^2=t^3$ is isomorphic to the subalgebra $A
:= K[x^2, x^3]$ of the polynomial algebra $K[x]$ (via $s\mapsto
x^3$, $t\mapsto x^2$). For a polynomial $a\in K[x]$, let $a':=
\frac{da}{dx}$ and $a'(0):=\frac{da}{dx}(0)$. Then
%\marginpar{Map0}
\begin{equation}\label{}
A = \{ a\in K[x]\, | \, a'(0)=0\}.
\end{equation}
The polynomial algebra $K[x]$ is a monoid with respect to the
composition $\circ$ of functions. It follows from the chain rule,
$(a\circ b)'= a'(b) b'$, that %\marginpar{Map1}
\begin{equation}\label{Map1}
K[x]\circ A \subseteq A \;\; {\rm and}\;\; A \circ (x) \subseteq A
\end{equation}
where $(x)$ is the ideal of the polynomial algebra $K[x]$
generated by the element $x$. In particular, $(A , \circ )$ is a
semigroup  but not a monoid. Indeed, suppose that $e$ is an
identity of $A$ then $\deg (a) = \deg (e \circ a) = \deg (e) \deg
(a)$ for all elements $a\in A$, and so $\deg (e) =1$. But the
semigroup $A$ contains no element of degree 1, a contradiction.

Note that $A \cap K[x]^* = \emptyset$. So, each element of $A$ is
not a unit of the monoid $(K[x], \circ )$.

The next lemma gives a necessary and sufficient condition for a
composition of two polynomials to be an element of $A$.

\begin{lemma}\label{2Sep7}%\marginpar{2Sep7}
Let $K$ be a field of characteristic zero and $a, b \in K[x]$.
Then $a\circ b\in A$ iff either $b\in A$ or $b\not\in A$ and the
value $b(0)$ of the polynomial $b(x)$ at $x=0$ is a root of the
derivative $\frac{da}{dx}$ of $a$.
\end{lemma}

{\it Proof}. $a\circ b\in A$ iff $0=(a\circ b)'(0)= a'(b(0))b'(0)$
iff either $b'(0)=0$ or, otherwise, $a'(b(0))=0$ iff either $b\in
A$ or, otherwise, $b(0)$ is a root of $a'$. $\Box $

$\noindent $

Let $\IrrA$ and $\IrrK[x]$ be the sets of irreducible elements of
the semi-groups $A$ and $K[x]$ respectively. The set $\IrrA$ is
the disjoint union of its two subsets $\CC$ and $\CD$ where
$$ \CC := \IrrA \cap \IrrK[x] = \{ p\in \IrrK[x] \, | \,
p'(0)=0\}$$ and $\CD := \IrrA \backslash \CC$. So, the set $\CC$
contains precisely all the irreducible elements of $K[x]$ that
belong to the semi-group $A$, and the set $\CD$ contains precisely
all the irreducible elements of $A$ which are {\em reducible} in
$K[x]$. Below, Proposition \ref{8Sep7} states a necessary and
sufficient condition for an irreducible element of $A$ to belong
to the set $\CC$ or $\CD$. First, let us give some definitions.

$\noindent $

 For a polynomial $a\in K[x]$,  let $\CR (a)$ and $\Dec (a)$ be, respectively,
  the set of its  roots
  and the set of all possible decompositions into irreducible factors in
  $K[x]$. For an element $a\in A$, let $\Dec_A(a)$ be the set of
  all possible decompositions into irreducible factors in $A$. If
  $p_1\circ \cdots \circ p_r\in \Dec (a)$ then
  $$a'=(p_1\circ \cdots \circ p_r)'=p_1'(p_2\circ \cdots \circ
  p_r) \cdot p_2'(p_3\circ \cdots \circ
  p_r)  \cdots p_{r-1}'(p_r)\cdot p_r', $$
and so %\marginpar{Rap}
\begin{equation}\label{Rap}
\CR (a') =\CR (p_1'(p_1\circ \cdots \circ
  p_{r-1}))\cup \cdots \cup \CR ( p_{r-1}'(p_r))\cup \CR (p_r').
\end{equation}
Let  $$\CE (a):= \cup_{p_1\circ \cdots \circ p_r\in \Dec (a)}\CR
(p_r').$$ By the very definition, the set  $\CE (a) $ is a subset
of $\CR (a')$. In particular, the set $\CE (a)$ is a finite set.
 In general, $\CE (a) \neq \CR (a')$. For each element $p\in
\IrrK[x]$, $ q\in \IrrA$ and $\l \in \CR(q')$, we have the
inclusions (where $K^*:=K\backslash \{ 0\}$)
$$ K[x]^*\circ p \circ K[x]^* \subseteq \IrrK[x]\;\; {\rm and}\;\;
K[x]^*\circ q \circ (\l + K^*x) \subseteq \IrrA .$$ In particular,
$K[x]^*\circ q \circ K^*x \subseteq \IrrA $ and $K[x]^*\circ q
\circ (\l +x) \subseteq \IrrA $.

\begin{proposition}\label{8Sep7}%\marginpar{8Sep7}
Let $p\in A\backslash K$. Then
\begin{enumerate}
\item $p\in \CC$ iff $p\in \IrrK[x]$ and $p'(0)=0$. \item  $p\in
\CD$ iff $p\not\in \CC$ and, for each decomposition $p_1\circ
\cdots \circ p_r\in \Dec (p)$, $(p_2\circ \cdots \circ
p_r)'(0)\neq 0$.
\end{enumerate}
\end{proposition}

{\it Proof}. 1. This is obvious.

2. $(\Rightarrow )$ Suppose that $p\in \CD$. Then, obviously,
$p\not\in \CC$. Suppose that $(p_2\circ \cdots \circ p_r)'(0)= 0$
for some decomposition $p_1\circ \cdots \circ p_r\in \Dec (p)$, we
seek a contradiction. Let $\l$ be a root of the polynomial $p_1'$.
The elements
$$ q_1:= p_1\circ (x+\l_1) \;\;{\rm and } \;\;  q_2:= (x-\l_1)^{-1}\circ p_2
\circ\cdots \circ p_r  $$  belong to the semi-group $A$, and
$$ p=q_1\circ q_2.$$
This contradicts to the irreducibility of the element $p$.
Therefore, $(p_2\circ \cdots \circ p_r)'(0)\neq 0$.

$(\Leftarrow )$ Suppose that $p\not\in \CC$ and, for each
decomposition $p_1\circ \cdots \circ p_r\in \Dec (p)$, $(p_2\circ
\cdots \circ p_r)'(0)\neq 0$. Suppose that the element $p$ is
reducible, i.e. $p=a\circ b$ for some elements $a,b\in A\backslash
K$, we seek a contradiction. Fix decompositions $p_1\circ \cdots
\circ p_s\in \Dec (a)$ and $p_{s+1} \circ \cdots \circ p_r\in \Dec
(b)$. Then $p= p_1\circ \cdots \circ p_r$ and  $(p_{s+1}\circ
\cdots \circ p_r)'(0)=0$ since $b\in A$, and so $(p_2\circ \cdots
\circ p_r)'(0)=0$ (by the chain rule), a contradiction. So, the
element $p$ is irreducible in $A$, hence $p\in \CD$ since
$p\not\in \CC$.  $\Box $

$\noindent $

The following two corollaries give a method of construction of
elements of the set $\CD$. In particular, they show that the set
$\CD$ is a non-empty set.
\begin{corollary}\label{a9Sep07}%\marginpar{a9Sep07}
Suppose that an element $q$ of $A$ is a composition $p_1\circ \cdots
\circ p_r$ of irreducible factors $p_i\in \IrrK[x]$ such that $r\geq
2$, $(p_2\circ \cdots \circ p_r)'(0)\neq 0$ and
$$\Dec (q)
=\{ (p_1\circ u_1)\circ (u_1^{-1}\circ p_2\circ u_2)\circ \cdots
\circ (u_{r-1}^{-1}\circ p_r)\, | u_1, \ldots , u_{r-1}\in
K[x]^*\}.$$ Then $q\in \CD$.
\end{corollary}

{\it Proof}. Since $r\geq 2$, $q\not\in \CC$. By the assumption,
for each decomposition $q_1\circ \cdots \circ q_r\in \Dec (q)$, we
can find elements $u_1, \ldots , u_{r-1}\in K[x]^*$ such that
$$ q_1=p_1\circ u_1, \;\; q_2=u_1^{-1}\circ p_2\circ u_2, \ldots , \;\;
q_r= u_{r-1}^{-1}\circ p_r.$$ Now, $\CR ((q_2\circ \cdots \circ
q_r)')= \CR ((u_1^{-1}\circ p_2\circ \cdots \circ p_r)')= \CR ((
p_2\circ \cdots \circ p_r)')\not\ni 0$. By Proposition
\ref{8Sep7}.(2), $q\in \CD$.  $\Box $

$\noindent $

Note that any sufficiently generic irreducible polynomials $p_1,
\ldots , p_r\in \IrrK[x]$ $(r\geq 2)$ with $p_1\circ \cdots \circ
p_r\in A$ satisfy the assumptions of Corollary \ref{a9Sep07}. For
example, take generic polynomials $p_1, \ldots , p_r\in K[x]$ such
that $(p_1\circ \cdots \circ p_r)'(0)=0$ and  $(p_2\circ \cdots
\circ p_r)'(0)\neq 0$ then all $p_i\in \IrrK[x]$ and $p_1\circ
\cdots \circ p_r\in \CD$.

\begin{corollary}\label{10Sep7}%\marginpar{10Sep7}
Let $r\geq 2$ be a natural number. For each natural  number $i=1,
\ldots , r$, let $p_i=\sum_{j=0}^{n_i} a_{ij} x^j\in K[x]$ be a
polynomial of prime degree $n_i\geq 5$. Suppose that $a_{11}:=-
\sum_{j=2}^{n_1}ja_{1j}(p_2\circ \cdots \circ p_r(0))^{j-1}$ and
that all the elements $a_{ij}$ of the field $K$ with  $(i,j) \neq
(1,1)$ are algebraically independent over the field of rational
numbers  $\Q$. Then $p_1\circ \cdots \circ p_r\in \CD$. In
particular, $\CD \neq \emptyset$.
\end{corollary}

{\it Proof}. The definition of the element $a_{11}$ means that
$p_1'((p_2\circ \cdots \circ p_r)(0))=0$. This implies that
$(p_1\circ\cdots\circ p_r)'(0)=0$, and so $p_1\circ\cdots\circ
p_r\in A$. Next, we show that the assumption of Corollary
\ref{a9Sep07} hold. The polynomials $p_i$ are irreducible since
their degrees are prime numbers. The elements $a_{ij}$, $i=2,
\ldots , r$, $j=1, \ldots , n_i$, are algebraically independent
over $\Q$, hence $(p_2\circ \cdots \circ p_r)'(0)\neq 0$. Suppose
that
$$\Dec
(p_1\circ \cdots \circ p_r)\neq \{ (p_1\circ u_1)\circ
(u_1^{-1}\circ p_2\circ u_2)\circ \cdots \circ (u_{r-1}^{-1}\circ
p_r)\, | u_1, \ldots , u_{r-1}\in K[x]^*\},$$
 we seek a contradiction. Then, by the second  theorem of Ritt-Levi, there
 exists a pair $(p_i, p_{i+1})$ and elements $\alpha , \beta , \g
 \in K[x]^*$ such that the pair $(\alpha \circ p_i\circ \beta ,
 \beta^{-1} \circ p_{i+1} \circ \g )$ is one of the three types:
\begin{eqnarray*}
&(a)&  (T_{n_i}, T_{n_{i+1}}),  \\
&(b)&   (x^{n_i} , x^rg(x^{n_i})),\;\;\;
 r+n_i\deg (g) = n_{i+1}, \\
 &(c)&  (x^r g^{n_{i+1}},
x^{n_{i+1}}), \;\;  r+n_{i+1} \deg (g) = n_i. \\
\end{eqnarray*}
For each polynomial $f\in K[x]$, let $C(f)$ be the subfield of $K$
generated by its  coefficients over $\Q$. In the case (a) (resp.
(b)) $p_i= \alpha^{-1} \circ T_{n_i}\circ \beta^{-1}$ (resp. $p_i=
\alpha^{-1} \circ x^{n_i}\circ \beta^{-1}$. On the one hand, the
transcendence degree $\trdeg \, C(p_i)= n_i\geq 5$, on the other
hand,  $\trdeg \, C(\alpha^{-1} \circ T_{n_i}\circ \beta^{-1})\leq
4$ (resp. $\trdeg \, C(\alpha^{-1} \circ x^{n_i}\circ
\beta^{-1})\leq 4$), a contradiction. Similarly, in the case (c),
$p_{i+1} = \beta \circ x^{n_{i+1}}\g^{-1}$, and so
$$ 5\leq \trdeg \, C(p_{i+1}) = \trdeg \, C(\beta \circ
x^{n_{i+1}}\g^{-1})\leq 4, $$ a contradiction. These
contradictions mean that the assumptions of Corollary
\ref{a9Sep07} hold for the element $p_1\circ \cdots \circ p_r$,
and so $p_1\circ \cdots \circ p_r\in \CD$. In particular, $\CD$ is
a non-empty set.  $\Box $

$\noindent $

The next lemma shows that, in general, the first theorem of J. F.
Ritt does not hold for the cusp.

\begin{lemma}\label{b9Sep7}%\marginpar{b9Sep7}
In general, the number of irreducible polynomials in decomposition
into irreducible polynomials  of an element of $A$ is non-unique.
Moreover, it can vary greatly.
\end{lemma}

{\it Proof}. Let $p\in \CD$ and $q\in \IrrA$. Consider their
composition $a:= p\circ q$.  Fix a decomposition $p_1\circ \cdots
\circ p_r\in \Dec (p)$, and then, for each $i=1, \ldots , r$, fix
a root, say $\l_i$, of the polynomial $p_i$. Consider the elements
of $\CC$:
$$a_1:= p_1\circ (x+\l_1), a_2:= (x-\l_1)^{-1}\circ p_2\circ (x+\l_2),
\ldots , a_r:= (x-\l_{r-1})^{-1}\circ p_r\circ (x+\l_r).$$ Then
$a_{r+1}:=(x-\l_r)^{-1}\circ q \in \IrrA$ and
$$a= p\circ q = a_1\circ \cdots \circ a_r\circ a_{r+1}$$
are two irreducible decompositions for the element $a$ with
distinct numbers of irreducible factors.  $\Box $

$\noindent $

Lemma \ref{b9Sep7} means that both theorems of J. F. Ritt fails
badly for the cusp. However,  we can describe a procedure of how
to obtain all irreducible decompositions of any given element of
$A$. Let $a\in A\backslash K$. Take any decomposition $p_1\circ
\cdots \circ p_r\in \Dec (a)$. Suppose that it is possible to
insert brackets
$$ (\ldots ) \circ (\ldots ) \circ \cdots \circ (\ldots)$$
in such a way that inside the brackets are irreducible elements of
$A$ (in principal, this can be  checked using Proposition
\ref{8Sep7}). It gives an irreducible decomposition for the element
$a$ in $A$. Moreover, all irreducible decompositions of the element
$a$ in $A$ can be obtained in this way.

$\noindent $

{\bf Proof of Theorem \ref{A15Sep7}.}

$\noindent $

We keep the notation of Theorem \ref{A15Sep7}. So, $a\in A\backslash
K$ with $l_A(a) = l(a)$, and $X,Y \in \Max (a)$. We have to show
that the decomposition $Y$ can be obtained from the decomposition
$X$ using some of the transformations (Adm), ($\CC a$), $\CC b $) or
$(\CC c$). We call these transformations the {\em cusp}
transformations. Note that $\Max ( a) \subseteq \Dec (a)$, and so
$X,Y \in \Dec (a)$. Let $X', Y'\in \Max (a)$. We write $X'\sim_A Y'$
if the decomposition $Y'$ can be obtained from the decomposition
$X'$ by using the cusp transformations. The relation $\sim_A$ on the
set $\Max (a)$ is an equivalence relation since the cusp
transformations are reversible. This means that the inverse of  a
transformation of the type (Adm) or ($\CC a$) is a transformation of
the type (Adm) or ($\CC a$) respectively; and the inverse of a
transformation of the type ($\CC b$) or ($\CC c$) is a
transformation of the type ($\CC b$) or ($\CC c$) respectively. We
write $X'\sim_\CC Y'$ if the decomposition $Y'$ is obtained from the
decomposition $X'$ by a single cusp transformation. Theorem
\ref{A15Sep7} means that the set $\Max (a)$ is an {\em equivalence
class} under the equivalence relation $\sim_A$, i.e. the equivalence
relation $\sim_A$ on $\Max (a)$ coincides with the equivalence
relation $\sim$, by the second theorem of Ritt-Levi (the equivalence
relation $\sim$ is defined in the proof of Theorem \ref{b31Aug7}).
We write $X'\sim_R Y'$ if $Y'$ is obtained from $X'$ by a single
Ritt transformation.

Let $r:= l_A(a) = l(a)$. Since $X,Y \in \Max (a)$, we have
$$ X= p_1\circ \cdots \circ p_r\;\; {\rm and }\;\; Y= q_1\circ \cdots \circ
q_r$$
 for some irreducible polynomials $p_i, q_i\in \CC$.

 {\em Case} $(\alpha )$: $K[x]^* p_r = K[x]^* q_r$, {\em i.e.}
 $q_r= \alpha \circ p_r$ {\em for some polynomial} $\alpha \in
 K[x]^*$. Let $b :=p_1\circ \cdots \circ p_{r-1}$. Then $b\circ p_r= a=q_1\circ \cdots \circ
q_r= q_1\circ \cdots \circ (q_{r-1}\circ \alpha ) \circ p_r$. By
(A4), we can delete $p_r$ at both ends of the chain of equalities
above, and the result is
$$ b =p_1\circ \cdots \circ p_{r-1}=q_1\circ \cdots \circ (q_{r-1}\circ \alpha
).$$ By Corollary \ref{c31Aug7}, the decomposition $V:= q_1\circ
\cdots \circ (q_{r-1}\circ \alpha )\in \Dec (b)$ can be obtained
from the decomposition $U:=p_1\circ \cdots \circ p_{r-1}\in \Dec
(b)$ by applying, say $t$, Ritt transformations
$$U=U_0\sim_RU_1\sim_RU_2\sim_R\cdots \sim_RU_t=V.$$
Then the decomposition $Y= V\circ p_r$ can be obtained from the
decomposition $X= U\circ p_r$ by applying cusp transformations of
the type (Adm) in the following way. First, we have the elements
of the set $\Dec (a)$:
$$ X=W_0:= U_0\circ p_r, \ldots , W_i:= U_i\circ p_r, \ldots ,
W_t:= U_t\circ p_r, \; W_{t+1}:=Y.$$ An important fact is that the
last element of all decompositions, that is $p_r$, is an element
of $A$. Let $U_i:= P_1\circ \cdots \circ P_{r-1}$ where $P_1,
\ldots , P_{r-1}\in \IrrK[x]$. Fort each polynomial $P_j$, fix a
$P_j$-admissible element, say $u_{ij}$, of $K[x]^*$, and consider
the decomposition
$$ W_i^*= P_i^*\circ \cdots \circ P_r^*\in \Max (a)$$
where
$$P_1^*:= P_1\circ u_{i1}, P_2^*:= u^{-1}_{i1}\circ P_2\circ u_{i2}, \ldots
, P_{r-1}^*:=u^{-1}_{i, r-2}\circ P_{r-1}\circ u_{i,r-1},
P_r^*:=u^{-1}_{i,r-1}\circ p_r.$$ It is obvious that the
decomposition $W_i^*$ is obtained from the decomposition $W_i$ by
applying $r-1$ transformations of the type (Adm). Let  $\Adm
(u_{i1}, \ldots , u_{i, r-1})$ denote their composition (in
arbitrary order since the transformations commute). We assume that
for $i=0,t+1$ all the $u$'s are equal to $x$. This means that the
transformation $\Adm (x, \ldots , x)$ is the identity
transformation, and, obviously, $W_0^* = W_0=X$ and $W_{t+1}^* =
W_{t+1}= Y$. So, there is the chain of elements of the set $\Max
(a)$:
$$ X=W_0^*,\;\;  W_1^*,\;\;  \ldots ,\;\;  W_t^*,\;\;  W_{t+1}^*=Y.$$
For each natural number $i=1, \ldots , t+1$, the decomposition
$W_i^*$ is obtained from the decomposition $W_{i-1}^*$ by applying
cusp transformations of the type (Adm):
$$ \Adm (u^{-1}_{i-1, 1}\circ u_{i1}, \ldots , u^{-1}_{i-1,r- 1}\circ u_{i,
r-1}).$$ Therefore, $X\sim_A Y$.

{\em Case} $(\beta )$: $K[x]^*p_r\neq K[x]^* q_r$. By Corollary
\ref{c31Aug7}, this means that $p_r=\l_{r-1}^{-1}\circ \pi \circ
\l_r$ for some units $\l_{r-1}, \l_r\in K[x]^*$ such that $\l_r$
is $\pi$-admissible and the polynomial $\pi$ is one of the
following types:

\begin{eqnarray*}
 &(a)& \pi =T_l, \;\; {\rm where} \;\; l \;\; {\rm is \; an \; odd \;  prime \; number}, \\
 &(b)& \pi =x^sg(x^p),  \;\; {\rm where} \;\; s\geq 1, \;\; g(x)\in K[x]\backslash K, \;\; g(0)\neq 0, \;\; p \;\; {\rm
 is\; a \;  prime \; number}, \\
 &(c)& \pi =x^p, \;\; {\rm where} \;\; p \;\; {\rm is \; a \;  prime \; number}.
\end{eqnarray*}

{\it Remark}. We exclude  the situation when $s=0$ in the case (b)
since otherwise we would have the case (c) due to irreducibility
of the element $\pi$ and the equality $g(x^p) = g(x) \circ x^p$.

We consider the three cases separately and label them respectively
as $(\beta a)$, $(\beta b )$ and $(\beta c)$.

{\em Case } $(\beta a)$: $\pi = T_l$ {\em where $l$ is an odd
prime number}. By the second theorem of Ritt-Levi, the element
$q_r$ in the decomposition $Y= q_1\circ \cdots \circ q_r$ must be
of the type $\mu \circ T_m\circ \l_r$ for some prime number $m$
such that $m\neq l$ (see Case $(\beta )$) where $\l_r$ is
necessarily a $T_m$-admissible polynomial and $\mu \in K[x]^*$. If
$\nu$ is the only root of the polynomial $\l_r$ then
$$\nu \in \CR (T_l')\cap \CR
(T_m')=\emptyset\;\;\;\; ({\rm Lemma} \;\ref{RTkTl}.(2)),$$
 a
contradiction. Therefore, this case is impossible.

{\em Case } $(\beta b)$: $\pi =x^sg(x^p)$  {\em (as in the case
(b) above)}. Then for the element $q_r$ there are two options
either $q_r\in K[x]^*\circ T_k\circ \l_r$ for some prime number
$k$ or, otherwise, $q_r\in K[x]^* \circ x^q\circ \l_r$ for some
prime number $q$. For $k\neq 2$, the first option is not possible
since by interchanging $X$ and $Y$ we would have the impossible
Case $(\beta a)$ (recall that the cusp transformations are
reversible). For $k=2$,  $T_2= (-1+2x)\circ x^2$, and so  we have,
in fact, only the second option, i.e. $q_r= \mu \circ x^q\circ
\l_r$ for some unit $\mu \in K[x]^*$. This means that the
invariant number
$$n_{\CP , q}\geq 1.$$
 Let $i$ be the {\em greatest } index
such that $p_i\in K[x]^*\circ x^q\circ K[x]^*$. In this case, we
call the element $p_i$  the {\em largest } $x^q$ in the
decomposition $X$ denoted $L(X)$. The decompositions
$$ H(X):= p_1\circ \cdots \circ p_{i-1} \;\; {\rm and }\;\; T(X):=
p_{i+1}\circ \cdots \circ p_r$$
 are called the {\em head}  and the {\em tail} of the
 decomposition $X$ respectively. The invariance of the number $n_{\CP ,
 q}$ means that we can control the largest $x^q$ under Ritt
 transformations. The largest $x^q$ remains unchanged under a Ritt
 transformation either of the head or the tail of $X$, and it
 moves to the right or left by one point if the largest $x^q$ is
 involved in the Ritt transformation of the type (b) or (c) from
 the Introduction
 respectively.

 Let $p_i= \l_{i-1}^{-1}\circ x^q\circ \l_i$ for some units
 $\l_{i-1}, \l_i\in K[x]^*$. Then the tail $T(X)$ of $X$ has clear
 structure. We claim that {\em there exist units} $\l_{i+1} , \ldots ,
 \l_{r-2}\in K[x]^*$ {\em such that}
 $$ p_j= \l_{j-1}^{-1} \circ \pi_j\circ \l_j, \;\; j=i+1, \ldots ,
 r-1, $$
{\em  where $\pi_j$ is either $x^n$ for a prime number $n$ or,
otherwise,
 $x^tf(x^q)$ for some $t\geq 1$ and  $f(x) \in K[x]$ such that  $\deg (f) \geq
 1$ and} $f(0)\neq 0$. The decomposition $Y$ is obtained from the
 decomposition $X$ by several Ritt transformations
 $$ X=X_0\sim_R X_1\sim_R\cdots \sim_R X_k\sim_R\cdots\sim_R X_m=Y.$$
Using the explicit form of Ritt transformations the claim  follows
easily by the backward induction on $k$ starting with the obvious
case $k=m-1$.

Using the claim we can produce $r-i$ cusp transformations
$$ X= Z_i\sim_CZ_{i+1}\sim_C\cdots \sim_C Z_r$$
such that on each step the largest $x^q$ moves one point to the
right, and the last irreducible element in the decomposition $Z_r$
is $q_r=\mu\circ x^q\circ \l_r$. On the first step,
$Z_i\sim_CZ_{i+1}$, the cusp transformation changes the triple
$$ (p_i, p_{i+1}, p_{i+2}) =(\l_{i-1}^{-1} \circ x^q\circ \l_i,\l_i^{-1} \circ \pi_{i+1}\circ \l_{i+1}, p_{i+2})$$
into the triple
$$
(p_i^*, p_{i+1}^*, p_{i+2}^*) =
\begin{cases}
(\l_{i-1}^{-1} \circ x^n, x^q, \l_{i+1}\circ p_{i+2})& \text{if $\pi_{i+1} = x^n$},\\
(\l_{i-1}^{-1} \circ [x^tf^q]\circ \nu , \nu \circ x^q, \l_{i+1}\circ p_{i+2})& \text{if  $\pi_{i+1} =x^tf(x^q)$},\\
\end{cases}$$
provided $i+1<r$ where $\nu \in K[x]^*$ is $x^tf^q$-admissible. If
$i+1=r$, the cusp transformation $Z_{r-1}\sim_CZ_r$ changes the
pair
$$ (p_{r-1}, p_r) =(\l_{r-2}^{-1} \circ x^n\circ \l_{r-1},\l_{r-1}^{-1} \circ [x^s h(x^q)]\circ \l_r)$$
into the pair
$$ (p_{r-1}^*, p_r^*) =(\l_{r-2}^{-1} \circ [x^sh^q],x^n\circ \l_r)$$
where $h(x^q) = g(x^p)$. The remaining cusp transformations are
defined by the same formulae as above by changing the index $i$
accordingly. Now, the decompositions $Z_r$ and $Y$ satisfy the
assumption of the case $(\alpha )$, and so $Z_r\sim_A Y$. Now,
$X\sim_A Z_r$ and $Z_r\sim_A Y$, and so  $X\sim_A Y$.

{\em Case } $(\beta c)$: $\pi =x^p$  {\em (as in the case (c)
above)}. The element $q_r$ has the form $\mu \circ
\widetilde{\pi}\circ \l_r$ where for the element $\widetilde{\pi}$
we have the same three options (a), (b) or (c) as for the element
$\pi$. Interchanging $X$ and $Y$, we reduce the cases (a) and (b)
for the element $\widetilde{\pi}$ to the cases (a) and (b) for
$\pi$ which have been considered already. For the last case,
$\widetilde{\pi}=x^q$, we repeat word for word the arguments of
the case $(\beta b)$ starting from the claim there. The proof of
Theorem \ref{A15Sep7} is complete. $\Box$

$\noindent $

{\bf Proof of Theorem \ref{15Sep7}.}

$\noindent $

Theorem \ref{15Sep7} follows easily from  the first theorem of J. F.
Ritt (or from Theorem \ref{A15Sep7} and the definition of the cusp
transformations, i.e. the transformations (Adm), $(\CC a)$, $(\CC
b)$ and $(\CC c)$). $\Box $

$\noindent $

\begin{proposition}\label{21Sep7}%\marginpar{21Sep7}
In general, Theorem \ref{15Sep7} does not hold for irregular
elements.
\end{proposition}

{\it Proof}. Let $m$ and $n$ be distinct prime numbers, $g(x)$ and
$h(x)$ be non-scalar polynomials of $K[x]$ such that $h(0)\neq 0$,
$k:= s+n \deg (g)$ and $l:= 1+m\deg (h)$ are prime numbers for
some natural number $s\geq 2$. Then the degrees of the polynomials
$x^n$, $x^sg(x^n)$ and $ xh(x^m)$ are prime numbers. Hence, the
polynomials $x^n$, $x^sg(x^n)$ and $x^sg^n$ are elements of the
set $\IrrA$, and $ xh(x^m)\in \IrrK[x]\backslash A$. It is obvious
that
$$ p:= [ x^sg(x^n)]\circ [ xh(x^m)],\; q:=x^n\circ [
xh(x^m)]\in \CD , $$ and the element  $a:= x^n\circ [
x^sg(x^n)]\circ [ xh(x^m)]\in A$ is irregular since $h(0)\neq 0$.
Then $$ a = x^n\circ p = x^sg^n\circ q \in \Dec_A(a),$$ $(\deg
(x^n) , \deg (p))= (n,kl)$ and $(\deg (x^sg^n), \deg (q))=
(k,nl)$. Since $k>n$, we have $(n,kl) \neq (k,nl)$ and  $(n,kl)
\neq (nl, k)$. This means that Theorem \ref{15Sep7} does not hold
for the irregular element $a$.  $\Box $

$\noindent $

In general, for an element $a$ of $A$ there exists a decomposition
$p_1\circ \cdots \circ p_t\in \Dec_A(a)$ with $t<l_A(a)$, i.e.
%\marginpar{MaDe}
\begin{equation}\label{MaDe}
\Max (a) \neq \Dec_A (a).
\end{equation}

{\it Example}. Let $k$ be an odd prime number, $g$ be a non-scalar
polynomial of $K[x]$ such that $l:= s+2 \deg (g)$ is a prime
number for some natural number $s\geq 2$. Let $\l$ be a root of
the trigonometric polynomial $T_k$. Consider the element $a:=
[x^sg^2]\circ T_k\circ T_2\in A$. The elements
$$ p_1:= x^sg^2, \;\; p_2:= T_k\circ (x+\l) \;\; {\rm and } \;\;
p_3:= (x-\l ) \circ T_2$$
 of the algebra $A$ are irreducible since their degrees are prime
 numbers. Let $q_1:= T_2$. Note that $q_2:= [ x^sg(x^2)]\circ
 T_k\in \CD$ since $T_k\in (x) \backslash (x^2)$ and $s\geq 2$.
 Then
 $$ a = p_1\circ p_2\circ p_3= q_1\circ q_2\in \Dec_A(a).\;\; \Box $$

$\noindent $

For an element $a$ of $A$, the number $\df (a) : = l(a) - l_A(a)$
is called the {\em defect} of the element $a$. The element $a$ is
irregular iff $\df (a) >0$. For each root $\l$ of the derivative
$a'$ of a polynomial $a$ of $K[x]$, the number
$$\ind_a(\l ) : = \max \{ i \, | \, \exists \; p_1\circ \cdots \circ
p_r\in \Dec (a) \; {\rm such \; that}\;\; p_i'(p_{i+1}\circ \cdots
\circ p_r\circ x)(0)=0\}$$ is called the {\em index} of $\l$. If
$a\in A$ then %\marginpar{lAind}
\begin{equation}\label{lAind}
l_A(a) = \ind_a(0).
\end{equation}
To prove this fact note that it is obvious that $l_A(a) \leq
\ind_a(0)$. For $i:= \ind_a(0)$, let us fix a decomposition
$p_1\circ \cdots \circ p_r\in \Dec (a)$ with $p_i'(p_{i+1}\circ
\cdots \circ p_r\circ x)(0)=0$. For each $j=1, \ldots , i-1$, let
$u_j$ be a $p_j$-admissible element of $K[x]^*$. The elements
$$ q_1:=p_1\circ u_1, \; q_2:= u_1^{-1} \circ p_2\circ u_2, \ldots
, q_{i-1}:= u_{i-2}^{-1} \circ p_{i-1} \circ u_{i-1} , \; q_i:=
u_{i-1}^{-1} \circ p_i\circ \cdots \circ p_r$$ belong to the
algebra $A$, and $a= q_1\circ \cdots \circ q_i$. Hence, $l_A(a)
\geq \ind_a(0)$. This establishes the equality (\ref{lAind}).

$\noindent $

For each element $a$ of $A$ with $i:= \ind_a(0)$, let
$$ \Dec (a, 0):= \{ p_1\circ \cdots \circ p_r\in \Dec (a) \, | \, p_i'(p_{i+1}\circ
\cdots \circ p_r\circ x)(0)=0\} .$$ The next lemma gives all the
 decompositions  of maximal length for each element  of $A$.

\begin{lemma}\label{22Sep7}%\marginpar{22Sep7}
Let $a$ be an element of $A$ and $i:= \ind_a(0)$. Then
\begin{eqnarray*}
\Max (a) &= & \{(p_1\circ u_1) \circ (u_1^{-1} \circ p_2\circ u_2)
\circ \cdots \circ
(u_{i-2}^{-1} \circ p_{i-1}\circ u_{i-1}) \circ  (u_{i-1}^{-1}\circ  p_i\circ\cdots \circ  p_r) \, | \,   \\
 & & p_1\circ \cdots \circ p_r\in \Dec (a, 0), \; u_j\in K[x]^*\;\; {\rm is } \; p_j-{\rm
 admissible}\}.
\end{eqnarray*}
\end{lemma}

{\it Proof}. It is obvious that the RHS $\subseteq \Max (a)$. On
the other hand, if $q_1\circ \cdots \circ q_i\in \Max (a)$ then
$q_1\circ \cdots \circ q_i\in $ the RHS. It suffices to put $p_j=
q_j$ and $u_j= x$. $\Box $

$\noindent $

By Lemma \ref{22Sep7}, if the element $a$ of $A$ is irregular and
$q_1\circ \cdots \circ q_i\in \Max (a)$ then necessarily $q_1,
\ldots , q_{i-1} \in \CC$ and $q_i\in \CD$.

$${\bf Acknowledgements}$$

The paper was finished during the author's visit to the IHES.
 Support and hospitality of the IHES is greatly acknowledged. The
 author would like to thank M. Zieve for comments and interesting
 discussions.

Department of Pure Mathematics

University of  Sheffield

Hicks Building

Sheffield S3~7RH

UK

email: v.bavula@sheffield.ac.uk

$\noindent $

IHES

Le Bois-Marie

35, Route de Chartes

F-91440 Bures-sur-Yvette

France

 email: bavula@ihes.fr

\end{document}